\documentclass[10pt,reqno]{amsart}
\usepackage{amsmath}
\usepackage{amssymb}
\usepackage{amsthm}

\newcommand{\sectie}[1]{\setcounter{equation}{0}\section{#1}}

\newcommand{\nabp}{\nab \hspace{-1.05ex}
\rule[.5ex]{.2ex}{.8ex}   \hspace{1.05ex}}

\newcommand{\golf}{\tilde{\rule{0ex}{1.7ex}}\,}

\newcommand{\Mhh}{\hat{M}
\hspace{-1.05ex}\hat{\rule{0ex}{2.0ex}}\hspace{1.05ex}}
\newcommand{\dehh}{\hat{\de}
\hspace{-.95ex}\hat{\rule{0ex}{2.05ex}}\hspace{.95ex}}
\newcommand{\vfihh}{\hat{\vfi}
\hspace{-.55ex}\hat{\rule{0ex}{1.45ex}}\hspace{.55ex}}
\newcommand{\lahh}{\hat{\la}
\hspace{-.8ex}\hat{\rule{0ex}{2.0ex}}\hspace{.8ex}}
\newcommand{\Gahh}{\hat{\Ga}
\hspace{-.8ex}\hat{\rule{0ex}{2.0ex}}\hspace{.8ex}}
\newcommand{\sdehh}{\hat{\sde}
\hspace{-.46ex}\hat{\rule{0ex}{1.95ex}}\hspace{.46ex}}

\newcommand{\spat}{\hspace{4ex}}

\newcommand{\flip}{\raisebox{.45ex}[0pt][0pt]{$\chi$}}

\newcommand{\Ga}{\Gamma}

\newcommand{\deh}{\hat{\Delta}}
\newcommand{\psih}{\hat{\psi}}
\newcommand{\vfih}{\hat{\vfi}}
\newcommand{\tauh}{\hat{\tau}}

\newcommand{\Rh}{\hat{R}}
\newcommand{\Sh}{\hat{S}}
\newcommand{\nabph}{\hat{\nabp}}

\newcommand{\lah}{\hat{\la}}
\newcommand{\Gah}{\hat{\Ga}}

\newcommand{\sih}{\hat{\si}}
\newcommand{\sdeh}{\hat{\delta}}

\newcommand{\nab}{\nabla}

\newcommand{\Jh}{\hat{J}}

\newcommand{\nabh}{\hat{\nab}}

\newcommand{\cI}{{\mathcal I}}

\newcommand{\cN}{{\mathcal N}}
\newcommand{\cM}{{\mathcal M}}

\newcommand{\ot}{\otimes}

\newcommand{\la}{\Lambda}

\newcommand{\om}{\omega}
\newcommand{\io}{\iota}
\newcommand{\vfi}{\varphi}

\newcommand{\ga}{\Gamma}
\newcommand{\sde}{\delta}
\newcommand{\de}{\Delta}

\newcommand{\si}{\sigma}
\newcommand{\Mfi}{{\mathcal M}_{\vfi}}
\newcommand{\Nfi}{{\mathcal N}_{\vfi}}

\newcommand{\Nps}{{\mathcal N}_{\psi}}
\newcommand{\Mpsi}{{\mathcal M}_{\psi}}
\newcommand{\Npsi}{{\mathcal N}_{\psi}}
\newcommand{\lafi}{\la_\vfi}
\newcommand{\laps}{\la_\psi}

\newcommand{\pifi}{\pi_\vfi}
\newcommand{\pips}{\pi_\psi}

\newcommand{\C}{\mathbb C}
\newcommand{\R}{\mathbb R}

\newcommand{\N}{\mathbb N}

\newcommand{\cT}{{\mathcal T}}

\newcommand{\dual}{M_*}
\newcommand{\duals}{M_*^\sharp}

\newcommand{\cst}{\text{C}$\hspace{0.1mm}^*$}

\newcommand{\Ext}{\hspace{-1.5ex}\raisebox{-0.5ex}[0pt][0pt]{\scriptsize\fontshape{n}\selectfont ext}}
\newcommand{\Mfih}{\mathcal{M}_{\vfih}}
\newcommand{\Op}{\raisebox{0.9ex}[0pt][0pt]{\scriptsize\fontshape{n}\selectfont op}\,}
\newcommand{\deop}{\de \hspace{-.3ex}\raisebox{0.9ex}[0pt][0pt]{\scriptsize\fontshape{n}\selectfont op}}
\newcommand{\dpr}{^{\prime\prime}}
\newcommand{\Si}{\Sigma}
\newcommand{\hoedje}{\,\hat{\rule{0ex}{1ex}}\;}

\newtheorem{definition}{Definition}[section]
\newtheorem{proposition}[definition]{Proposition}
\newtheorem{lemma}[definition]{Lemma}
\newtheorem{corollary}[definition]{Corollary}

\newtheorem{theorem}[definition]{Theorem}

\newtheorem{notation}[definition]{Notation}

\newcommand{\Mr}{M_{\text{\tiny c}}}
\newcommand{\der}{\Delta_{\text{\tiny c}}}
\newcommand{\vfir}{\vfi_{\text{\tiny c}}}
\newcommand{\psir}{\psi_{\text{\tiny c}}}
\newcommand{\taur}{\tau^{\text{\tiny c}}}
\newcommand{\sir}{\si^{\text{\tiny c}}}
\newcommand{\Nfir}{{\mathcal N}_{\vfir}}
\newcommand{\Mfir}{{\mathcal M}_{\vfir}}

\newcommand{\lar}{\la_{\text{\tiny c}}}

\newcommand{\Rr}{R_{\text{\tiny c}}}
\newcommand{\Sr}{S_{\text{\tiny c}}}

\newcommand{\csot}{\ot_{\text{\tiny c}}}

\newcommand{\tekst}[1]{\;\,\text{#1}\;\,}
\newcommand{\Nfih}{\mathcal{N}_{\hat{\varphi}}}

\begin{document}

\begin{center}

\begin{center}
\Huge\bf Locally quantum groups in the von Neumann algebraic setting  \end{center}

\bigskip\bigskip

\rm Johan Kustermans \& Stefaan Vaes

\medskip

Department of Mathematics

KU Leuven

Celestijnenlaan 200B

B--3001 Heverlee

Belgium

\medskip

e-mails : johan.kustermans@wis.kuleuven.ac.be \& stefaan.vaes@wis.kuleuven.ac.be

\bigskip\bigskip

\bf May 2000 \rm
\end{center}

\bigskip\medskip

\subsection*{Abstract}
In this paper we complete in several aspects the picture of locally compact
quantum groups. First of all we give a definition of a locally compact quantum
group in the von Neumann algebraic setting and show how to deduce from it a
C$^*$-algebraic quantum group. Further we prove several results about locally
compact quantum groups which are important for applications, but were not yet
settled in our paper \cite{KV1}. We prove a serious strengthening of the left
invariance of the Haar weight, and we give several formulas connecting the
locally compact quantum group with its dual. Loosely speaking we show how the
antipode of the locally compact quantum group determines the modular group and
modular conjugation of the dual locally compact quantum group.

\bigskip

\section*{Introduction}

\medskip

Building on the work of Kac \& Vainerman \cite{VK}, Enock \& Schwartz
\cite{ES}, Baaj \& Skandalis \cite{BS}, Woronowicz \cite{Wor5} and Van Daele
\cite{VD} a precise definition of a locally compact quantum group was recently
introduced by the authors in \cite{KV1}, see \cite{KV2} and \cite{KV3} for an
overview. For an overview of the historic development of the theory we refer to
\cite{KV3} and the introduction to \cite{KV1}. Because commutative
C$^*$-algebras are always of the form $C_0(X)$, where $X$ is a locally compact
space and $C_0(X)$ denotes the C$^*$-algebra of continuous functions on $X$
vanishing at infinity, arbitrary C$^*$-algebras are sometimes thought of as the
algebra of continuous functions vanishing at infinity on a (non-existing)
locally compact quantum space. For this reason the C$^*$-algebra framework is
the most natural one to study locally compact quantum groups. The most general
commutative example of a locally compact quantum group is $C_0(G)$ with
comultiplication $\de: C_0(G) \rightarrow C_b(G \times G)$ given by $(\de
f)(x,y) = f(xy)$, where $G$ is a locally compact group and $C_b$ denotes the
algebra of continuous bounded functions. This philosophy is followed in
\cite{KV1} where we defined \lq reduced C$^*$-algebraic quantum groups\rq\ as
the proper notion of a locally compact quantum group in the C$^*$-algebra
framework. As already explained, the most general commutative example is
$C_0(G)$ where $G$ is a locally compact group. Further the theory unifies
compact quantum groups and Kac algebras and it includes known examples as the
quantum Heisenberg group, quantum $E(2)$-group, quantum Lorentz group and
quantum $az+b$-group. Within this theory one can construct a dual reduced
C$^*$-algebraic quantum group and prove a Pontryagin duality theorem.

\smallskip

On a technical level it is often more easy to work with von Neumann algebras
rather than C$^*$-algebras, certainly when dealing with weights. So, already in
\cite{KV1}, we associated with every reduced C$^*$-algebraic quantum group a
von Neumann algebraic quantum group and we used it to prove several results on
the C$^*$-algebra level. The first aim of this paper is to give an intrinsic
definition of a von Neumann algebraic quantum group and to associate with it,
in a canonical way, a reduced C$^*$-algebraic quantum group. This can be
thought of as the quantum analogue of the classical result of Weil (see
\cite[Appendice I]{Weil}), stating that every group with an invariant measure
has a unique topology turning it into a locally compact group.

\smallskip

A second aim of this paper is to prove some new results on both C$^*$-algebraic
and von Neumann algebraic quantum groups, which are indispensable for
applications. In our definition of either C$^*$-algebraic or von Neumann
algebraic quantum groups we assume the existence of left and right invariant
weights. But the property of invariance we assume is quite weak, and in this
paper we show how a much stronger notion of invariance can be proved. The same
kind of result is stated for Kac algebras in \cite{EScros}, but not proved. The
first proof was given by Zsid\'o in \cite{Zsido} (see also remark 18.23 in
\cite{Stra}). This stronger invariance property is needed whenever an action of
a von Neumann algebraic quantum group on a von Neumann algebra appears: see
\cite{Kust} and \cite{SV}, but also \cite{EScros} for Kac algebra actions, and
it will certainly be useful in future investigations as well.

\smallskip

Further we will complete the picture of the quantum group and its dual with
several formulas giving a link between the antipode of the quantum group and
the modular theory of its dual. Roughly speaking we obtain that $$\hat{T}^*
\la(x) = \la \bigl( S(x^*) \bigr)$$ for nice $x \in M$, where $M$ is the von
Neumann algebraic quantum group, $\la$ is the GNS-map of the left invariant
weight $\vfi$ on $M$, $S$ is the antipode and $\hat{T}$ is the operator
appearing in the modular theory of the left invariant weight $\vfih$ on the
dual von Neumann algebraic quantum group: it is the closure of $\lah(\om)
\mapsto \lah(\om^*)$ where $\lah$ is the GNS-map of $\vfih$. To these results
and formulas will be referred in further research, see e.g. \cite{Kust} and
\cite{SV}.

\bigskip

We end this introduction with some conventions and references concerning
weights and operator valued weights.

\medskip

We assume that the reader is familiar with the theory of normal semi-finite
faithful weights (in short, n.s.f. weights) on von Neumann algebras.
Nevertheless, let us fix some notations. So let $\vfi$ be a n.s.f. weight on a
von Neumann algebra $M$. Then we define the following sets:
\begin{enumerate}
\item $\cM_\vfi^+ = \{ \, x \in M^+ \mid \vfi(x) < \infty \, \}$, so $\cM_\vfi^+$ is a hereditary cone in $M^+$,
\item $\cN_\vfi = \{ \, x \in M \mid x^*x \in \cM_\vfi^+ \, \}$, so $\cN_\vfi$ is a left ideal in $M$,
\item $\cM_\vfi =$ the linear span of $\cM_\vfi^+$ in $M$, so $\cM_\vfi$ is a $^*$-subalgebra of $M$.
\end{enumerate}
There exists a unique linear map $F : \cM_\vfi \rightarrow \C$ such that $F(x)
= \vfi(x)$ for all $x \in \cM_\vfi^+$. For all $x \in \cM_\vfi$, we set
$\vfi(x) = F(x)$.

\smallskip

A GNS-construction for $\vfi$ is a triple $(H_\vfi,\pifi,\lafi)$, where
$H_\vfi$ is a Hilbert space, $\pifi : M \rightarrow B(H_\vfi)$ is a normal
$^*$-homomorphism and $\lafi : \cN_\vfi \rightarrow H_\vfi$ is a
$\si$-strong$^*$ closed linear map with dense range such that (1) $\langle
\lafi(x), \lafi(y) \rangle = \vfi(y^* x)$ for all $x,y \in \cN_\vfi$ and (2)
$\lafi(x\, y) = \pifi(x) \lafi(y)$ for all $x \in M$ and $y \in \cN_\vfi$. As
usual we introduce the closed densely defined linear operator $T$ in $H_\vfi$
as the closure of the map $\lafi(x) \mapsto \lafi(x^*)$ for $x \in \Nfi \cap
\Nfi^*$. Making the polar decomposition $T=J\nab^{\frac{1}{2}}$ of $T$ we
obtain the modular operator $\nab$ and modular conjugation $J$ of $\vfi$ with
respect to the GNS-construction $(H_\vfi,\pifi,\lafi)$.

\medskip

Consider two von Neumann algebras $M$, $N$. Let $\vfi$ be a n.s.f. weight on
$M$ with GNS-construction $(H_\vfi,\pifi,\lafi)$ and let $\psi$ be a n.s.f.
weight on $N$ with GNS-construction $(H_\psi,\pips,\laps)$. The tensor product
weight $\vfi \ot \psi$ is a n.s.f. weight on $M \ot N$ (see e.g. definition 8.2
of \cite{Stra} for a definition). This tensor product weight has a
GNS-construction $(H_\vfi \ot H_\psi,\pifi \ot \pips,\lafi \ot \laps)$ where
$\lafi \ot \laps : \cN_{\vfi \ot\psi} \rightarrow H_\vfi \ot H_\psi$ is the
$\si$-strong$^*$ closure of the algebraic tensor product $\lafi \odot \laps :
\cN_{\vfi} \odot \cN_{\psi} \rightarrow H_\vfi \ot H_\psi$.

\medskip

Let $M$ be any von Neumann algebra. For the definition of the extended positive
part $M^+\Ext$ we refer to definition 1.1 of \cite{H4}. For $T \in M^+\Ext$ and
$\om \in M_*^+$, we set $\langle T , \om \rangle  = T(\om) \in [0,\infty]$.
Recall that there exists an embedding $M^+  \hookrightarrow M^+\Ext : x \mapsto
x^\sharp$ such that $\langle x^\sharp , \om \rangle = \om(x)$ for all $x \in
M^+$ and $\om \in M_*^+$. We will use this embedding to identify $M^+$ as a
subset of $M^+\Ext$.

\smallskip

Consider a von Neumann algebra $M$ and a von Neumann subalgebra $N$ of $M$. The
definition of an operator valued weight from $M$ to $N$ is given in definition
2.1 of \cite{H4}.

\medskip

Now consider two von Neumann algebras $M$ and $N$ and a n.s.f. weight $\vfi$ on
$M$.  We identify $N$ with $\C \ot N$ as a von Neumann subalgebra of $M \ot N$
to get into the framework of operator valued weights. The operator valued
weight $\vfi \ot \io : (M \ot N)^+ \rightarrow N^+\Ext$ is defined in such a
way that for $x \in (M \ot N)^+$, we have that $$\om\bigl((\vfi \ot
\io)(x)\bigr) = \vfi\bigl((\io \ot \om)(x)\bigr)\ .$$ As for weights we define
the following sets:
\begin{enumerate}
\item $\cM_{\vfi \ot \io}^+ = \{ \, x \in (M \ot N)^+ \mid (\vfi \ot \io)(x) \in N^+ \, \}$, so $\cM_{\vfi \ot \io}^+$ is a hereditary cone of $(M \ot N)^+$,
\item $\cN_{\vfi \ot \io} = \{ \, x \in M \ot N \mid x^*x \in \cM_{\vfi \ot \io}^+ \, \}$, so $\cN_{\vfi \ot \io}$ is a left ideal in $M \ot N$,
\item $\cM_{\vfi \ot \io} =$ the linear span of $\cM_{\vfi \ot \io}^+$ in $M \ot N$, so $\cM_{\vfi \ot \io}$ is a $^*$-subalgebra of $M \ot N$.
\end{enumerate}
There exists a unique linear map $G : \cM_{\vfi \ot \io} \rightarrow N$ such
that $G(x) = (\vfi \ot \io)(x)$ for all $x \in \cM_{\vfi \ot \io}^+$. For all
$x \in \cM_{\vfi \ot \io}$, we set $(\vfi \ot \io)(x) = G(x)$. Let $a \in
\cM_\vfi$ and $b \in N$. Then it is easy to see that $a \ot b$ belongs to
$\cM_{\vfi \ot \io}$ and $(\vfi \ot \io)(a \ot b) = \vfi(a) \, b$.

\smallskip

Thanks to the remark after lemma 1.4 of \cite{H4}, we also have the following
characterization of $\cM_{\vfi \ot \io}^+$: Let $x \in (M \ot N)^+$, then $x$
belongs to $\cM_{\vfi \ot \io}^+$ $\Leftrightarrow$ $\vfi((\io \ot \om)(x)) <
\infty$ for all $\om \in M^+_*$.

Let $x \in \cN_{\vfi \ot \io}$ and $\om \in N_*$. The inequality $(\io \ot
\om)(x)^* (\io \ot \om)(x) \leq \|\om\| \, (\io \ot |\om|)(x^* x)$ will imply
that $(\io \ot \om)(x) \in \cN_\vfi$ and $$\|\lafi((\io \ot \om)(x))\| \leq
\|\om\|\,\| (\vfi \ot \io)(x^*x)  \|^{\frac{1}{2}} \ .$$

\smallskip

When $L$ is some set of elements of a space we denote by $\langle L \rangle$
the linear span of $L$ and by $[L]$ the closed linear span. The symbol $\ot$
will denote either a von Neumann algebraic tensor product or a tensor product
of Hilbert spaces and $\io$ will denote the identity map. Finally we use the
symbol $\flip$ to denote the flip map from $M \ot N$ to $N \ot M$, where $N$
and $M$ are von Neumann algebras. We use $\Si$ to denote the flip map from $H
\ot K$ to $K \ot H$ when $H$ and $K$ are Hilbert spaces.

\medskip

\bigskip

\sectie{Von Neumann algebraic quantum groups} \label{sectie1}

\medskip

We state the definition of a von Neumann algebraic quantum group and discuss
how the \cst-algebraic theory can be translated to the von Neumann algebraic
setting. The major difference between both approaches is the absence of density
conditions in the definition of von Neumann algebraic quantum groups: these
will follow automaticly!

\smallskip

\begin{definition}
Consider a von Neumann algebra $M$ together with a unital normal
$^*$-homomorphism $\de : M \rightarrow M \ot M$ such that $(\de \ot \io)\de =
(\io \ot \de)\de$. Assume moreover the existence of
\begin{enumerate}
\item  a n.s.f. weight $\vfi$ on $M$ that is left invariant:
$\vfi((\om \ot \io) \de(x)) = \vfi(x) \om(1)$ for all $\om \in M_*^+$ and $x
\in \Mfi^+$.
\item  a n.s.f. weight $\psi$ on $M$ that is right invariant:
$\psi((\io \ot \om) \de(x)) = \psi(x) \om(1)$ for all $\om \in M_*^+$ and $x
\in \Mpsi^+$.
\end{enumerate}
Then we call the pair $(M,\de)$ a von Neumann algebraic quantum group.
\end{definition}

\smallskip

In the next part of this section we will list the essential properties of these
quantum groups. Most of the time the proofs in \cite{KV1} can be easily
translated to the von Neumann algebraic setting by replacing the norm and
strict topology in the considerations by the $\si$-strong$^*$ topology.
However, some care has to be taken to prove the density conditions and we will
discuss this in detail.

\smallskip

For the rest of this section we fix a von Neumann algebraic quantum group
$(M,\de)$. Without loss of generality, we may and  will assume that $M$ is in
standard form with respect to a Hilbert space $H$.

At the same time we fix a n.s.f. left invariant weight $\vfi$ on $(M,\de)$
together with a GNS-construction $(H,\io,\la)$ (which is possible because $M$
is in standard form). We let $\nab$ denote the modular operator and $J$ the
modular conjugation of $\vfi$ with respect to this GNS-construction
$(H,\io,\la)$.

\smallskip

We also choose a n.s.f. right invariant weight $\psi$ on $(M,\de)$ together
with a GNS-construction $(H,\io,\Gamma)$ (later on, we will introduce some
canonical choice for $\psi$ and $\Gamma$).

\medskip

By left invariance of $\vfi$, we get  that $(\om \ot \io)\de(x) \in \cN_\vfi$
and $\|\la((\om \ot \io)\de(x))\| \leq \|\om\|\,\|\la(x)\|$ for all $x \in
\cN_\vfi$ and $\om \in M_*$. Arguing as in result 2.6 of \cite{KV1}, we also
get  for $a,b \in \cN_\psi$ and  $x \in \Nfi$ that $(\psi \ot \io)(\de(b^* x)(a
\ot 1)) \in \Nfi$ and $\|\la\bigl((\psi \ot \io)(\de(b^* x)(a \ot 1))\bigr)\|
\leq \|\ga(a)\|\,\|\ga(b)\|\, \|\la(c)\|$.

\smallskip

Along the way to the proof of theorem \ref{def.thm1}, one also translates
proposition 3.15 of \cite{KV1}, giving rise to the important equalities
\begin{eqnarray}
H & = & [\,\la((\om \ot \io)\de(x)) \mid x \in \cN_\vfi, \om \in M_*\,]
\label{eq1}
\\ & = & [\,\la\bigl((\psi \ot \io)(\de(b^* x)(a \ot 1))\bigr) \mid x \in \cN_\vfi, a,b \in \Npsi\,] \ . \label{eq2}
\end{eqnarray}

\medskip

The left invariance of $\vfi$ implies that $\de(y)(x \ot 1) \in \cN_{\vfi \ot
\vfi}$ for all $x,y \in \cN_\vfi$ and that $$\langle (\la \ot \la)(\de(y_1)(x_1
\ot 1)) , (\la \ot \la)(\de(y_2)(x_2 \ot 1)) \rangle = \langle \la(x_1) \ot
\la(y_1) , \la(x_2) \ot \la(y_2) \rangle$$ for all $x_1,x_2,y_1,y_2 \in \Nfi$.
The proof of the next result is an easy translation of the proof of theorem
3.16 of \cite{KV1}.

\begin{theorem} \label{def.thm1}
There exists a unique unitary element $W \in B(H \ot H)$ such that $W^* (\la(x)
\ot \la(y)) = (\la \ot \la)(\de(y)(x \ot 1))$ for all $x,y \in \cN_\vfi$.
\end{theorem}

\medskip

It should be noted that $(\om \ot \io)(W^*) \la(x) = \la((\om \ot \io)\de(x))$
for all $x \in \Nfi$ and $\om \in B(H)_*$ (see e.g. result 2.10 of \cite{KV1}).
Using the commutant theorem for the tensor product of von Neumann algebras,
this implies that $W$ is a unitary element in $M \ot B(H)$.

\smallskip

Using the formula for $W^*$ above, one sees that $\de(x) = W^* (1 \ot x) W$ for
all $x \in M$. Applying the techniques of proposition 3.18 of \cite{KV1}, one
checks that $W$ satisfies the pentagonal equation: $W_{12} W_{13} W_{23} =
W_{23} W_{12}$. We call $W$ the multiplicative unitary of $(M,\de)$ with
respect to the GNS-construction $(H,\io,\la)$.

\medskip

It goes without saying that all these results also have their right invariant
counterparts. For later purposes we introduce the unitary element $V \in B(H)
\ot M$ such that $V(\Gamma(x) \ot \Gamma(y)) = (\Gamma \ot \Gamma)(\de(x)(1 \ot
y))$ for all $x,y \in \Npsi$. As in result 2.10 of \cite{KV1}, one proves that
\begin{equation}\label{sliceV}
(\om_{\ga(a),\ga(b)} \ot \io)(V^*) = (\psi \ot \io)(\de(b^*)(a \ot 1))
\quad\text{for all}\quad a,b \in \Npsi.
\end{equation}

\medskip

The proof of proposition 3.22 of \cite{KV1} survives the translation to the von
Neumann algebra setting. Combining this with equation (\ref{eq2}) we arrive at
the following conclusion.

\begin{proposition} \label{toep.antipode}
There exists a unique densely defined closed antilinear operator $G$ in $H$
such that $$\langle\,  \la\bigl((\psi \ot \io)(\de(x^*)(y \ot 1))\bigr) \mid
x,y \in \Nfi^*\,\Npsi\,  \rangle$$ is a core for $G$ and $$G \la\bigl((\psi \ot
\io)(\de(x^*)(y \ot 1))\bigr) = \la\bigl((\psi \ot \io)(\de(y^*)(x \ot
1))\bigr)$$ for $x,y \in \Nfi^*\,\Npsi$. We have moreover that $G$ is
involutive.
\end{proposition}

\medskip

By taking the polar decomposition of $G$, we get the following essential
operators in $H$.

\begin{notation} \label{notatieN}
We define $N = G^* G$, so $N$ is a strictly positive operator in $H$. We also
define the anti-unitary operator $I$ on $H$ such that $G = I \,
N^{\frac{1}{2}}$.
\end{notation}

Because $G$ is involutive, we have that $I = I^*$, $I^2 = 1$ and $I \, N \, I =
N^{-1}$.

\medskip

A careful analysis of the proof of proposition 5.5 of \cite{KV1} reveals that
this result remains true in the present setting. By equation~(\ref{sliceV}) and
the techniques used in the proof of proposition 5.8 of \cite{KV1}, this is
equivalent to saying that $$(\om_{v,w} \ot \io)(V^*) \, G \subseteq G \,
(\om_{w,v} \ot \io)(V^*) \hspace{1cm} \text{and} \hspace{1cm} (\om_{v,w} \ot
\io)(V) \, G^*  \subseteq G^* \, (\om_{w,v} \ot \io)(V)$$ for all $v,w \in H$.
Hence, appealing to the proof of result 5.10 of \cite{KV1}, we arrive at the
vital commutation relation
\begin{equation}
V (\nab_\psi \ot N) = (\nab_\psi \ot N) V \ , \label{eq3}
\end{equation}
where $\nab_\psi$ denotes the modular operator of $\psi$ with respect to the
GNS-construction $(H,\io,\ga)$.

\medskip

Up till now, we did not need the density conditions that are present in the
definition of reduced \cst-algebraic quantum groups (see definition 4.1 of
\cite{KV1}). This is the case because we were only working on the Hilbert space
level for which the relevant density conditions are already established in
equations (\ref{eq1}) and (\ref{eq2}). In order to further develop the theory
along the lines of \cite{KV1}, we will now prove similar density conditions on
the level of the von Neumann algebra $M$. The idea of the proof is taken from
\cite[2.7.6]{ES}.

\begin{proposition} \label{def.prop1}
Denoting by $^-$ the $\si$-strong$^*$ closure we have
\begin{align*}
M &= \langle (\om \ot \io)\de(x) \mid x \in M, \om \in M_* \rangle^- \\ &=
\langle (\io \ot \om)\de(x) \mid x \in M, \om \in M_* \rangle^-
\\ &= \bigl\{ (\om \ot \io)(V) \mid \om \in B(H)_* \bigr\}^- \ .
\end{align*}
\end{proposition}

\begin{proof}
Define $\cT_\psi$ to be the Tomita $^*$-algebra of $\psi$. From formula
(\ref{sliceV}) it follows that
\begin{align*}
\bigl\{ (\om \ot \io)(V) \mid \om \in B(H)_* \bigr\}^- &= \langle (\psi \ot
\io) \bigl( (ca^* \ot 1) \de(b) \bigr) \mid a,b \in \Npsi, c \in \cT_\psi
\rangle^- \\ &= \langle (\psi \ot \io) \bigl( (a^* \ot 1) \de(b)
(\si_{-i}^\psi(c) \ot 1) \bigr) \mid a,b \in \Npsi, c \in \cT_\psi \rangle^- \\
&= \langle (\om \ot \io)\de(x) \mid x \in M, \om \in M_* \rangle^- \ .
\end{align*}
Now we define $$M_r = \langle (\om \ot \io)\de(x) \mid x \in M, \om \in M_*
\rangle^- \ .$$ Because $V$ is a multiplicative unitary the linear space
$\bigl\{ (\om \ot \io)(V) \mid \om \in B(H)_* \bigr\}$ is an algebra that acts
non-degenerately on $H$. Because $M_r$ is clearly self-adjoint, we get that
$M_r$ is a von Neumann subalgebra of $M$. Working with the von Neumann
algebraic quantum group $(M,\flip \de)$ instead of $(M,\de)$ we obtain that
also $$M_l = \langle (\io \ot \om)\de(x) \mid x \in M, \om \in M_* \rangle^-$$
is a von Neumann subalgebra of $M$. Observe that it follows from the commutant
theorem for the tensor product of von Neumann algebras that $\de(x) \in M_l \ot
M_r$ for all $x \in M$.

\smallskip

Then we conclude from equation (\ref{eq3}) that it is possible to define a
one-parameter group $(\tau_t)_{t \in \R}$ of automorphisms of $M_r$ by
$\tau_t(x)=N^{-it} x N^{it}$ for all $x \in M_r$ and $t \in \R$. It also
follows from equation (\ref{eq3}) and the fact $\de(x) = V(x \ot 1)V^*$ for all
$x \in M$, that we have $\de(\si^\psi_t(x)) = (\si^\psi_t \ot \tau_{-t})\de(x)$
for all $x \in M$ and $t \in \R$, which makes sense because $\de(x) \in M \ot
M_r$. For the same reason we can write $$M_l = \bigl\{(\io \ot \om)\de(x) \mid
x \in M, \om \in (M_r)_* \bigr\}^-$$ and because $\si^\psi_t \bigl((\io \ot
\om)\de(x) \bigr) = (\io \ot \om \tau_t)\de(\si^\psi_t(x))$ for all $\om \in
(M_r)_*$ and $x \in M$, we get $\si_t^\psi(M_l) = M_l$ for all $t \in \R$. By
the right invariance of $\psi$ it follows that the restriction $\psi_l$ of
$\psi$ to $M_l$ is semifinite. By Takesaki's theorem (see e.g.
\cite[10.1]{Stra}) there exists a unique normal faithful conditional
expectation $E$ from $M$ to $M_l$ satisfying $\psi(x) = \psi_l(E(x))$ for all
$x \in M^+$. From \cite[10.2]{Stra} it follows that $E(x)P = PxP$ for all $x
\in M$, where $P$ denotes the orthogonal projection onto the closure of
$\Gamma(\Npsi \cap M_l)$. So the range of $P$ contains $\Gamma \bigl((\io \ot
\om)\de(x) \bigr)$ for all $\om \in M_*$ and $x \in \Npsi$. By the right
invariant version of equation (\ref{eq1}) we get that $P=1$. So $E$ is the
identity map and $M_l=M$.

\smallskip

Working with the von Neumann algebraic quantum group $(M,\flip \de)$ we obtain
$M=M_r$. We already proved that $M_r$ is the $\si$-strong$^*$ closure of
$\{(\om \ot \io)(V) \mid \om \in B(H)_* \}$ and so this concludes the proof of
the proposition.
\end{proof}

Because we have proved that $(M,\de)$ satisfies the above density conditions,
it is straightforward to translate the rest of the the proofs in \cite{KV1} to
the von Neumann algebraic setting. In the following part of this section, we
collect the most important results (we will not stick to the order as they
appear in \cite{KV1}).

\medskip

\subsection*{Uniqueness of invariant weights}

An essential result is the uniqueness of left and right invariant weights. If
$\theta$ is a normal semi-finite left invariant weight on $(M,\de)$, then there
exists a non-negative number $r$ such that $\theta = r \, \vfi$. A similar
result holds for right invariant weights.

\medskip

\subsection*{The antipode and its polar decomposition}

The antipode of our quantum group is defined through its polar decomposition.
There exists a strongly continuous one-parameter group $\tau$ on $M$  such that
$\tau_t(x) = N^{-it} x N^{it}$ for all $x \in M$ and $t \in \R$. At the same
time, we have a $^*$-anti-automorphism $R$ on $M$ such that $R(x) = I x^* I$
for all $x \in M$.

Then $R^2 = \io$, $R$ and $\tau$ commute and we define $S = R
\tau_{-\frac{i}{2}} = \tau_{-\frac{i}{2}} R$. Note that these 3 properties
determine the pair $R$, $\tau$ completely in terms of the map $S$. The map $S :
D(S) \subseteq M \rightarrow M$ is a $\si$-strongly$^*$ closed map with
$\si$-strong$^*$ dense domain and range that is determined by $(M,\de)$ through
the following so-called strong left invariance properties.

We have for all $a,b \in \Nfi$ that $(\io \ot \vfi)(\de(a^*)(1 \ot b)) \in
D(S)$ and
\begin{equation}
S\bigl((\io \ot \vfi)(\de(a^*)(1 \ot b))\bigr) = (\io \ot \vfi)((1 \ot
a^*)\de(b)) \ . \label{strong}
\end{equation}
The  space $\langle \, (\io \ot \vfi)(\de(a^*)(1 \ot b)) \mid a,b \in
\Nfi\,\rangle$ is a core for $S$. A similar result holds for right invariant
weights (see proposition 5.24 of \cite{KV1}). For other characterizations of
$S$ we refer to proposition 5.33 and corollary 5.34 of \cite{KV1}.

\medskip

We refer to $S$ as the antipode of the quantum group $(M,\de)$. The
one-parameter group $\tau$ is called the scaling group of $(M,\de)$, the map
$R$ is called the unitary antipode of $(M,\de)$.

There also exists a unique strictly positive number  $\nu$ such that $\vfi
\tau_t = \nu^{-t}\,\vfi$ for all $t \in \R$. We call $\nu$ the scaling constant
of $(M,\de)$. In connection with this relative invariance, it is useful to
define the strictly positive operator $P$ in $H$ such that $P^{it} \, \la(x) =
\nu^{\frac{t}{2}} \, \la(\tau_t(x))$ for all $t \in \R$ and $x \in \Nfi$. We
observe that $\tau_t(x) = P^{it} x P^{-it}$ for all $t \in \R$ and $x \in M$.

\medskip

\subsection*{The right Haar weight and the modular element}

Because we have the equation $\flip(R \ot R)\de = \de R$, we get that $\vfi R$
is a right invariant n.s.f. weight on $(M,\de)$. From now on we suppose that
$\psi = \vfi R$. Let $\si'$ denote the modular group of $\psi$. Remember that
$\si_t' = R \si_{-t} R$ for all $t \in R$. We have that $\vfi \si_t' =
\nu^t\,\vfi$, $\psi \si_t = \nu^{-t} \, \psi$ and $\psi \tau_t =
\nu^{-t}\,\psi$ for all $t \in \R$.

\smallskip

By the Radon Nikodym theorem 5.5 of \cite{SV2}, we get the existence of a
unique strictly positive element $\sde$ affiliated to $M$ such that
$\si_t(\sde) = \nu^t \,\sde$ for all $t \in \R$ and $\psi = \vfi_\sde$ (see
definition 1.5 of \cite{SV2} for the precise definition of $\vfi_\sde$).
Formally we have $\psi(x)=\vfi(\sde^{1/2} x \sde^{1/2})$. The element $\sde$ is
called the modular element of $(M,\de)$.

We have that $\de(\sde) = \sde \ot \sde$, $R(\sde) = \sde^{-1}$ and
$\tau_t(\sde) = \sde$ for all $t \in \R$.

\smallskip

Now we choose the GNS-construction $(H,\io,\Ga)$ for $\psi$ such that $\Ga =
\la_\sde$ (see the remarks before proposition 1.15 in \cite{KV1} for a precise
definition of $\la_\sde$). We denote the modular operator  of $\psi$ in this
GNS-construction by $\nabp$. Recall that $\nu^{\frac{i}{4}} \, J$ is the
modular conjugation of $\psi$ with respect to this same GNS-construction.

\medskip

\subsection*{The fundamental commutation relations}

A full-fledged theory of quantum groups would be impossible without the
following list of commutation relations.
\begin{enumerate}
\item The one-parameter groups $\tau$, $\si$ and $\si'$ commute pairwise.
\item For all $t \in \R$ we have
\begin{equation} \label{eq4}
\begin{array}{rclcrcl} \de \, \si_t  & =  & (\tau_t
\ot \si_t)\de  & \hspace{2cm}  & \de\, \si'_t & = & (\si'_t \ot \tau_{-t}) \de
\\ \de \,\tau_t & =  &(\tau_t \ot
\tau_t) \de & \hspace{2cm}  & \de \, \tau_t & = & (\si_t \ot \si'_{-t}) \de
\end{array}
\end{equation}
\item On the Hilbert space level, we get that
\begin{eqnarray}
(I \ot J) W  & = & W^* (I \ot J) \label{eq5}
\\ (N^{-1} \ot \nab) W & = & W (N^{-1} \ot \nab) \label{eq6}
\end{eqnarray}
\end{enumerate}

\medskip

\subsection*{The dual von Neumann algebraic quantum group}

The multiplicative unitary $W$ is manageable in the sense of \cite{Wor5} with
$P$ as the managing positive operator (see proposition 6.10 of \cite{KV1}). We
follow more or less chapter~3 of \cite{ES} to obtain the Haar weight on the
dual von Neumann algebraic quantum group (see section~8 of \cite{KV1} for the
proofs of the next results).

\begin{definition}
Define $\hat{M}$ to be the $\si$-strong$^*$ closure of the algebra $\{\,(\om
\ot \io)(W) \mid \om \in B(H)_*\,\}$. Then $\hat{M}$ is a von Neumann algebra
and there exists a unique unital normal $^*$-homomorphism $\deh : \hat{M}
\rightarrow \hat{M} \ot \hat{M}$ such that $\deh(x) = \Sigma W (x \ot 1) W^*
\Sigma$ for all $x \in \hat{M}$. The pair $(\hat{M},\deh)$ is again a von
Neumann algebraic quantum group, referred to as the dual of $(M,\de)$.
\end{definition}

\smallskip

The predual $M_*$ is a Banach algebra if we define the product such that $\om
\theta = (\om \ot \theta)\de$ for all $\om,\theta \in M_*$ (of course, $M_*$
should be thought of as the space of $L^1$-functions of $M$). Moreover, the map
$\lambda : M_* \rightarrow \hat{M} : \om \mapsto (\om \ot \io)(W)$ is an
injective morphism of algebras.

Let us recall the construction of the dual weight $\vfih$. First of all, we
define $$ \cI = \{\,\om \in M_* \mid \exists M \in \R^+ : |\om(x^*)| \leq M \,
\|\la(x)\| \text{ for all } x \in \Nfi\,\} \ .$$ By the Riesz theorem for
Hilbert spaces there exists for every $\om \in \cI$ a unique element $\xi(\om)
\in H$ such that $\om(x^*) = \langle \xi(\om) , \la(x) \rangle$ for all $x \in
\Nfi$.
 Then $\cI$ is a left ideal in $M_*$, the map $\cI \rightarrow  H : \om
\mapsto \xi(\om)$ is linear and $\lambda(\eta) \xi(\om) = \xi(\eta \om)$ for
all $\eta \in M_*$ and $\om \in \cI$.

\smallskip

There exists a unique $\si$-strong$^*$--norm closed linear map $\lah$, with
$\si$-strong$^*$ dense domain $D(\lah) \subseteq \hat{M}$, into $H$ such that
$\lambda(\cI)$ is a $\si$-strong$^*$--norm core for $\lah$ and
$\lah(\lambda(\om)) = \xi(\om)$ for all $\om \in \cI$. The dual weight $\vfih$
is the unique n.s.f. weight on $\hat{M}$ having the triple $(H,\io,\lah)$ as a
GNS-construction. It turns out that $\vfih$ is left invariant with respect to
$(\hat{M},\deh)$. We denote the modular group of $\vfih$ by $\sih$.

\smallskip

We denote the antipode, unitary antipode and scaling group of $(\hat{M},\deh)$
by $\Sh$, $\Rh$ and $\tauh$ respectively. The scaling constant of
$(\hat{M},\deh)$ is equal to $\nu^{-1}$. Define the right invariant n.s.f.
weight $\psih$ on $(\hat{M},\deh)$ as $\psih = \vfih \Rh$. The modular group of
$\psih$ will be denoted by $\sih'$. Denote the modular element of
$(\hat{M},\deh)$ by $\sdeh$. Referring to the fact that $\psih =
\vfih_{\sdeh}$, we define the GNS-construction $(H,\io,\hat{\Gamma})$ such that
$\hat{\Gamma} = \lah_{\sdeh}$.

\smallskip

The modular operator and modular conjugation of $\vfih$ with respect to
$(H,\io,\lah)$ will be denoted by $\nabh$ and $\hat{J}$ respectively. We will
denote the modular operator of $\psih$ with respect to $(H,\io,\hat{\Ga})$ by
$\nabph$.

It is also worth mentioning that $P^{it} \lah(x) = \nu^{-\frac{t}{2}} \,
\lah(\tauh_t(x))$ for all $t \in \R$ and $x \in \cN_{\vfih}$, which means in a
sense that $\hat{P}=P$.

\smallskip

Finally we mention that $M \cap \hat{M} = \C$.

\medskip

\subsection*{Pontryagin duality}

As in the previous paragraph, we can also construct the dual $(\Mhh,\dehh)$ of
$(\hat{M},\deh)$. Notice that the construction of the dual depends on the
choice of the GNS-construction of the left Haar weight. If we use the
GNS-construction $(H,\io,\lah)$ for the construction of the dual $(\Mhh,\dehh)$
, the Pontryagin duality theorem tells us that $(\Mhh,\dehh) = (M,\de)$. We
even have that $\vfihh = \vfi$ and $\lahh = \la$.

Since $(\Mhh,\dehh) = (M,\de)$, we get that  $\sdehh = \sde$. Hence $\Gahh =
\la_\sde = \Ga$.

\medskip

\subsection*{From von Neumann algebraic to \cst-algebraic quantum
groups}

In \cite{KV1}, we associated to any reduced \cst-algebraic quantum group a von
Neumann algebraic quantum group by taking the $\si$-strong$^*$ closure of the
underlying \cst-algebra in the GNS-space of a left Haar weight. In the last
part of this section we go the other way around by introducing a \cst-algebraic
quantum group.

To distinguish between von Neumann algebraic and \cst-algebraic tensor products
we will denote the minimal \cst-tensor product by $\csot$.

\begin{proposition}
Define $\Mr$ to be the norm closure of the space $\{\,(\io \ot \om)(W) \mid \om
\in B(H)_*\,\}$ and $\der$ to be the restriction of $\de$ to $\Mr$. Then the
pair $(\Mr,\der)$ is a reduced \cst-algebraic quantum group.
\end{proposition}
\begin{proof}
Because $W$ is manageable and $\der(x) = W^* (1 \ot x) W$ for all $x \in \Mr$,
propositions 1.5 and 5.1 of \cite{Wor5} imply that $\Mr$ is a \cst-algebra,
$\der$ is a non-degenerate $^*$-homomorphism from $\Mr$ into the multiplier
algebra of $\Mr \csot \Mr$, such that $(\der \csot \io)\der = (\io \csot
\der)\der$ and both $\der(\Mr)(\Mr \ot 1)$ and $\der(\Mr)(1 \ot \Mr)$ are dense
in $\Mr \csot \Mr$.

Now define $\vfir$ and $\psir$ to be the restriction of $\vfi$ and $\psi$ to
$\Mr^+$ respectively, giving you two faithful lower semi-continuous weights on
$\Mr$.

By equation (\ref{eq5}) we get that $(I \ot J)W(I \ot J) = W^*$, implying that
$R((\io \ot \om_{v,w})(W)) = (\io \ot \om_{J w, J v})(W)$ for all $v,w \in H$.
It follows that $R(\Mr) = \Mr$. Define $\Rr$ to be the restriction of $R$ to
$\Mr$. Then $\Rr$ is a $^*$-anti-automorphism of $\Mr$ satisfying $\flip(\Rr
\csot \Rr)\der = \der \Rr$. It is also clear that $\psir = \vfir \Rr$.

For $a,b \in \Npsi$ and $c \in \Nfi$, we have that $$(\psi \ot \io)(\de(b^*
c)(a \ot 1)) = R\bigl((\io \ot \vfi) ((1 \ot R(a))\de(R(c) R(b)^*))\bigr) =
R\bigl((\io \ot \om_{\la(R(c) R(b)^*),\la(R(a)^*)})(W^*)\bigr) \ ,$$ which
implies that $\Mr = [\,(\psi \ot \io)(\de(b^* c)(a \ot 1)) \mid a,b \in \Nps, c
\in \Nfi\,]$.

We know that $(\psi \ot \io)(\de(b^* c)(a \ot 1)) \in \Nfi$ and thus $(\psi \ot
\io)(\de(b^* c)(a \ot 1)) \in \Nfir$ for all $a,b \in \Npsi$ and $c \in \Nfi$.
It follows that $\vfir$ is densely defined.

Define $\lar$ to be the restriction of $\la$ to $\Nfir$. Equation (\ref{eq2})
guarantees that $\lar(\Nfir)$ is dense in $H$. Therefore $(H,\io,\lar)$ is a
GNS-construction for $\vfir$.

Equation (\ref{eq6}) tells us that $(N^{-1} \ot \nab) W = W (N^{-1} \ot \nab)$
implying that
\begin{equation} \label{ster1}
\tau_t((\io \ot \om_{v,w})(W)) = (\io \ot \om_{\nab^{it} v , \nab^{it} w})(W)
\end{equation}
for all $v,w \in H$ and $t \in \R$. Hence $\tau_t(\Mr) = \Mr$ for all $t \in
\R$. Define the one-parameter group $\taur$ on $\Mr$ by setting $\taur_t =
\tau_t\!\restriction_{\Mr}$ for all $t \in \R$. Notice that equation
(\ref{ster1}) implies that $\taur$ is norm continuous.

Since $\der(\Mr)(1 \ot \Mr)$ is a dense subset of $\Mr \csot \Mr$, we get that
$\Mr = [\,(\io \ot \om)\de(x) \mid \om \in B(H)_*,x \in \Mr\,]$. Equation
(\ref{eq4}) implies for all $t \in \R$, $\om \in B(H)_*$ and $x \in \Mr$ that
\begin{equation} \label{ster2}
\si_t((\io \ot \om)\de(x)) = (\io \ot \om \si_{t}')\de(\taur_t(x)).
\end{equation}
Therefore $\si_t(\Mr) = \Mr$ for all $t \in \R$ and we can define a one
parameter group $\sir$ on $\Mr$ by setting  $\sir_t =
\si_t\!\restriction_{\Mr}$ for all $t \in \R$. Equation (\ref{ster2}) implies
that $\sir$ is norm continuous.

By now it is clear that $\vfir$ is a KMS-weight on $\Mr$ (in the \cst-algebraic
sense) with $\sir$ as its modular group. Because $\psi = \vfi R$, we also get
that $\psir$ is a KMS-weight on $\Mr$.

\medskip

Take $\om \in (\Mr)^*_+$ and $x \in \Mfir^+$. Choose $\eta \in B(H)_*^+$. On
the \cst-algebra $\Mr$ we can make a GNS-construction for the positive
functional $\om$. This way we obtain a Hilbert space $K$, a non-degenerate
representation $\pi$ of $\Mr$ on $K$ and a (cyclic) vector $v \in K$ such that
$\om = \om_{v,v} \pi$. By theorem 1.5 of \cite{Wor5} we know that $W$ belongs
to the multiplier algebra of $\Mr \ot B_0(H)$, where $B_0(H)$ denotes the
\cst-algebra of compact operators on $H$. Hence the unitary $U$ defined by
$U:=(\pi \csot \io)(W)$ belongs to $B(K) \ot B(H)$. Define $\theta \in
B(H)_*^+$ by setting $\theta(x) = (\om_{v,v} \ot \eta)(U^* (1 \ot x) U)$ for
all $x \in B(H)$. Then
\begin{align*}
(\eta \ot \io)\de\bigl((\om \csot \io)(\der(x))\bigr) &= (\eta \ot
\io)\bigl((\om \csot \io \csot \io)((\der \csot \io)\der (x))\bigr) \\ &= (\eta
\ot \io)\bigl( (\om_{v,v} \ot \io \ot \io)(U^*_{12} \de(x)_{23} U_{12})) =
(\theta \ot \io)\de(x) \ .
\end{align*}
Therefore the left invariance of $\vfi$ implies that $(\eta \ot
\io)\de\bigl((\om \csot \io)(\der(x))\bigr)$ belongs to $\Mfi^+$ and
\begin{equation}\label{bijnaleft}
\vfi\bigl((\eta \ot \io)\de\bigl((\om \csot \io)(\der(x))\bigr)\bigr) =
\theta(1)\, \vfi(x) = \om(1) \, \eta(1) \, \vfir(x) \ .
\end{equation}
Translating proposition 5.15 of \cite{KV1} to the von Neumann algebra setting,
we now conclude that $(\om \csot \io)(\der(x))$ belongs to $\Mfi^+$ and
therefore to $\Mfir^+$.

Taking $\eta \in B(H)_*$ such that $\eta(1)=1$, equation (\ref{bijnaleft}) and
the left invariance of $\vfi$ imply  that $$\vfir\bigl( (\om \csot
\io)(\der(x))\bigr) = \vfi\bigl( (\om \csot \io)(\der(x))\bigr) =
\vfi\bigl((\eta \ot \io)\de\bigl((\om \csot \io)(\der(x))\bigr)\bigr) =  \om(1)
\, \vfir(x) \ .$$

So we have proven that $\vfir$ is left invariant in the sense of definition 2.2
of \cite{KV1}. Because $\flip(\Rr \csot \Rr)\de = \de \Rr$ and $\psir = \vfir
\Rr$ we also get that $\psir$ is right invariant. From all this we conclude
that $(\Mr,\der)$ is a reduced \cst-algebraic quantum group.
\end{proof}

\medskip

The GNS-construction $(H,\io,\lar)$ for $\vfir$ was obtained by letting $\lar$
be the restriction of $\la$ to $\Nfir$. By the definitions introduced at the
end of section 4 in \cite{KV1}, it is clear that this implies that $W$ is the
multiplicative unitary  of $(\Mr,\der)$ in this GNS-construction
$(H,\io,\lar)$.

Since $\sir_t$ and $\taur_t$ are restrictions of $\si_t$ and $\tau_t$
respectively, it is clear that $(\taur_t \ot \sir_t)\de = \de \sir_t$ for all
$t \in \R$, and so the density conditions imply that $\taur$ is the scaling
group of $(\Mr,\der)$. It follows that $\nu$ is also the scaling constant of
$(\Mr,\der)$. Letting $\Sr$ denote the antipode of $(\Mr,\der)$, proposition
5.33 of \cite{KV1} and its von Neumann algebraic counterpart imply that $\Sr
\subseteq S$. Since $\taur$ is the scaling group of $(\Mr,\der)$ and $\Rr$ was
obtained by restricting $R$ to $\Mr$, this implies that $\Rr$ is the unitary
antipode of $(\Mr,\der)$.

\smallskip

Following \cite{KV1}, we associate to the reduced \cst-algebraic quantum group
$(\Mr,\der)$ the von Neumman algebraic quantum group $(\tilde{M}_{\text{\tiny
c}},\tilde{\de}_{\text{\tiny c}})$ by letting $\tilde{M}_{\text{\tiny c}}$ be
the $\si$-strong$^*$ closure of $\Mr$ and defining $\tilde{\de}_{\text{\tiny
c}}$ to be the unique normal $^*$-homomorphism from $\tilde{M}_{\text{\tiny
c}}$ to $\tilde{M}_{\text{\tiny c}} \ot \tilde{M}_{\text{\tiny c}}$ extending
$\der$. It follows from proposition \ref{def.prop1} that
$(\tilde{M}_{\text{\tiny c}},\tilde{\de}_{\text{\tiny c}}) = (M,\de)$. We get
similar results for the extensions of the Haar weights, their modular groups,
the scaling group, the unitary antipode and the antipode itself.

\bigskip

\sectie{Commutation relations and related matters}

\medskip

In this section we  establish some useful technical properties about von
Neumann algebraic quantum groups that are often used when working in the
operator algebra approach to quantum groups. We start of by implementing the
scaling groups and unitary antipodes. Then we prove some results concerning the
dual and end by formulating some commutation relations.

\smallskip

We still have fixed a von Neumann algebraic quantum group $(M,\de)$ and we use
the notations introduced in section \ref{sectie1}.

\medskip

\begin{proposition} \label{implement}
The following properties hold $$\begin{array}{lclclcl} \tau_t(x) & = &
\nabh^{it} x \,\nabh^{-it} \text{\ \  for all } t \in \R, x \in M &
\hspace{1.5cm} & R(x) & = & \hat{J} x^* \hat{J} \text{\ \ for all } x \in M
\\ \tauh_t(x) & = & \nab^{it} x \,\nab^{-it} \text{\ \ for all } t \in \R, x \in \hat{M} & \hspace{1.5cm} &
\Rh(x) & = & J x^* J \text{\ \ for all } x \in \hat{M}
\end{array}$$
\end{proposition}
\begin{proof}
By propositions 8.17 and 8.25 of \cite{KV1}, we know that $\hat{R}(x) = J x^*
J$ for all $x \in \hat{M}$. Therefore the Pontryagin duality theorem guarantees
that also $R(x) = \hat{J} x^* \hat{J}$ for all $x \in M$.

Choose $x \in M$. By lemma 8.8 and proposition 8.9 of \cite{KV1} we know that
$\nabh^{it} = P^{it} J \sde^{it} J$, and so we get that $\nabh^{it} x\,
\nabh^{-it} = P^{it} J \sde^{it} J\, x \,J \sde^{-it} J P^{-it}$. But
Tomita-Takesaki theory tells us that $J \sde^{it} J$ belongs to $M'$, implying
that $\nabh^{it} x\, \nabh^{-it}  = P^{it} x\, P^{-it} = \tau_t(x)$.

Pontryagin duality allows us to conclude that $\tauh^{it}(x) = \nab^{it} x
\nab^{-it}$ for all $t \in \R$ and $x \in \hat{M}$.
\end{proof}

\medskip

We have that $(R \ot \hat{R})(W) = W$ and $(\tau_t \ot \tauh_t)(W) = W$ for all
$t \in \R$ (see the remarks before propositions 8.18 and 8.25 of \cite{KV1}).
Hence the next result.

\begin{corollary} \label{JenW}
We have the following commutation relations: $$W (\nabh \ot \nab) =  (\nabh \ot
\nab) W \hspace{2cm} \text{and} \hspace{2cm} W (\hat{J} \ot J) = (\hat{J} \ot
J) W^* \ .$$
\end{corollary}

\smallskip\smallskip

Notice that for the same reasons, $W (P \ot \nab) =  (P \ot \nab) W$ and $W
(\nabh \ot P) = (\nabh \ot P) W$.

\bigskip\medskip

In the next part, we  complete the picture of the dual. For this reason, let us
introduce a natural $^*$-algebra inside $\dual$.

\smallskip

\begin{definition}
Define the subspace $\duals$ of $\dual$ as $$\duals = \{\, \om \in \dual \mid
\exists \,\theta \in \dual : \theta(x) = \overline{\om}(S(x)) \text{ for all }
x \in D(S) \, \}\ .$$ We define the antilinear mapping $.^* : \duals
\rightarrow \duals$ such that $\om^*(x) = \overline{\om}(S(x))$  for all $\om
\in \duals$ and $x \in D(S)$. Then $\duals$ is a subalgebra of $\dual$ and
becomes a $^*$-algebra under the operation $.^*$\ .
\end{definition}

\smallskip

If $x \in D(S)$, then $S(x)^* \in D(S)$ and $S(S(x)^*)^* = x$ (which follows
from the corresponding property for $\tau_{-\frac{i}{2}}\,$). This implies that
$.^*$ is an involution on $\duals$. If $\om$ and $\eta$ are elements in $M_*$
such that $\om S$ and $\eta S$ are bounded and their unique continuous linear
extensions  $(\om S)\golf$ and $(\eta S)\golf$ belong to $M_*$, then $(\om
\eta)S$ is bounded and  $(\eta S)\golf\,(\om S)\golf$ is its unique continuous
linear extension (cfr lemma 5.25 of \cite{KV1}). From this, it follows that
$\duals$ is a subalgebra of $\dual$  and that $.^*$ is antimultiplicative.
Notice that, since $S$ can be unbounded, $\duals$ can be strictly smaller than
$\dual$.

\medskip\smallskip

For $\om \in B(H)_*$, the element $(\io \ot \om)(W)$ belongs to $D(S)$ and
$S((\io \ot \om)(W)) = (\io \ot \om)(W^*)$. The space $\{\,(\io \ot \om)(W)
\mid \om \in B(H)_*\,\}$ is moreover a $\si$-strong$^*$ core for $S$. (see
proposition 8.3 of \cite{KV1}). We use this characterization of $S$ to prove
the next result.

\begin{proposition}
The following holds :
\begin{trivlist}
\item[\ \,1.] $\duals = \{ \, \om \in \dual \mid \exists \,\theta \in \dual : \lambda(\om)^* = \lambda(\theta) \, \}$.
\item[\ \,2.] $\lambda(\om)^* = \lambda(\om^*)$ for all $\om \in \duals$.
\end{trivlist}
\end{proposition}
\begin{proof}
Take $\om \in \dual$. Then we have for all $\eta \in B(H)_*$ that
\begin{equation}\label{ster3}
\overline{\om}\bigl(S((\io \ot \eta)(W))\bigr)
 =  \overline{\om}((\io \ot \eta)(W^*))
= \eta((\om \ot \io)(W)^*) = \eta(\lambda(\om)^*) \ .
\end{equation}

\medskip

If $\om \in \duals$ then the formula above implies for all $\eta \in B(H)_*$
that $$\eta(\lambda(\om^*)) = \eta((\om^* \ot \io)(W)) = \om^*((\io \ot
\eta)(W)) = \eta(\lambda(\om)^*) \ , $$ and hence $\lambda(\om^*) =
\lambda(\om)^*$.

\medskip

Now suppose  that there exists $\theta \in \dual$ such that $\lambda(\om)^* =
\lambda(\theta)$. By formula (\ref{ster3}) above, we get for all $\eta \in
B(H)_*$ that $$\overline{\om}\bigl(S((\io \ot \eta)(W))\bigr) = \eta((\theta
\ot \io)(W)) = \theta((\io \ot \eta)(W)) \ . $$ Because such elements $(\io \ot
\eta)(W)$ form a $\si$-strong$^*$ core for $S$, we get $\overline{\om}(S(x)) =
\theta(x)$ for all $x \in D(S)$. So $\om$ belongs to $\duals$.
\end{proof}

\medskip

It is easy to prove that $\duals$ is dense in $\dual$ implying that the
$^*$-algebra $\lambda(\duals)$ is $\si$-strong$^*$ dense in $\hat{M}$. For
later purposes, we will need a result which gives a little bit more
information.

\begin{lemma}
The spaces $\cI \cap \duals$ and $(\cI \cap \duals)^*$ are dense in $\dual$ and
$\lambda(\cI \cap \duals)$ is a $\si$-strong$^*$--norm core for $\lah$.
\end{lemma}
\begin{proof}
Consider $\om \in \cI$. For every $n \in \N$ and $z \in \C$, we define
$\om(n,z) \in \dual$ as $$\om(n,z) = \frac{n}{\sqrt{\pi}} \int \exp(-n^2
(t+z)^2)\, \om \tau_t\, dt \ .$$ So we have for $x \in D(S)$ that $x \in
D(\tau_{-\frac{i}{2}})$ and thus
\begin{eqnarray*}
\overline{\om(n,z)}(S(x)) & = & \frac{n}{\sqrt{\pi}} \int \exp(-n^2
(t+\bar{z})^2)\, \overline{\om}\bigl(\tau_t(S(x))\bigr)\, dt  \\ & = &
\frac{n}{\sqrt{\pi}} \int \exp(-n^2(t +
\bar{z})^2)\,\overline{\om}\bigl(R(\tau_{t-\frac{i}{2}}(x))\bigr)\,dt \\ & = &
\frac{n}{\sqrt{\pi}} \int \exp(-n^2(t + \frac{i}{2} +
\bar{z})^2)\,\overline{\om}\bigl(R(\tau_t(x))\bigr)\,dt \ ,
\end{eqnarray*}
from which we conclude that $\om(n,z) \in \duals$ and
\begin{equation} \label{duals.eq1}
\om(n,z)^* = \frac{n}{\sqrt{\pi}} \int \exp(-n^2 (t+ \frac{i}{2}+\bar{z})^2)\,
\overline{\om} R \tau_t\, dt \ .
\end{equation}
It is easy to check that for every $t \in \R$ we have $\om \tau_t \in \cI$ and
$\xi(\om \tau_t) = \nu^{-\frac{t}{2}} \, P^{-it} \xi(\om)$. Therefore the
closedness of the mapping $\eta \mapsto \xi(\eta)$ implies that $\om(n,z) \in
\cI$ and
\begin{equation} \label{duals.eq2}
\xi(\om(n,z)) = \frac{n}{\sqrt{\pi}} \int \exp(-n^2 (t+z)^2)\,
\nu^{-\frac{t}{2}}\, P^{-it} \xi(\om)\, dt \ .
\end{equation}
\begin{trivlist}
\item[\ \,(1)] Let $\om \in \cI$. Then we have for every $n \in \N$ that $\om(n,0) \in \cI \cap \duals$. Clearly, $(\om(n,0))_{n=1}^\infty$ converges to  $\om$. Equation (\ref{duals.eq2}) implies that  $\bigl(\xi(\om(n,0))\bigr)_{n=1}^\infty$ converges to $\xi(\om)$. In other words, $\bigl(\lah(\lambda(\om(n,0)))\bigr)_{n=1}^\infty$ converges to $\lah(\lambda(\om))$.
Since $\cI$ is dense in $\dual$ and $\lambda(\cI)$ is a core for $\lah$, we
conclude that $\cI \cap \duals$ is dense in $\dual$ and that $\lambda(\cI \cap
\duals)$ is a $\si$-strong$^*$--norm core for $\lah$.
\item[\ \,(2)] Let $\om \in \cI$. Then we have for every $n \in \N$ that $\om(n,\frac{i}{2}) \in \cI \cap \duals$ and
$$\om(n,\frac{i}{2})^* = \frac{n}{\sqrt{\pi}} \int \exp(-n^2 t^2)\,
\overline{\om} R \tau_t\, dt $$ by equation (\ref{duals.eq1}). So
$(\om(n,\frac{i}{2})^*)_{n=1}^\infty$ converges to $\overline{\om}R$.  From
this all, we conclude that $(\cI \cap \duals)^*$ is dense in $\dual$.
\end{trivlist}
\end{proof}

\bigskip

\begin{proposition} \label{dual.prop20}
Define $\cI^\sharp = \{\,x \in \cI \cap \duals \mid x^* \in \cI \, \}$. Then
$\cI^\sharp$ is a $^*$-subalgebra of $\duals$ such that $\cI^\sharp$ is dense
in $\dual$ and $\lambda(\cI^\sharp)$ is a $\si$-strong$^*$--norm core for
$\lah$.
\end{proposition}
\begin{proof}
It is clear that $\cI^\sharp$ is a $^*$-subalgebra of $\duals$. Because $\cI$
is a left ideal in $\dual$, we get that $(\cI \cap \duals)^*(\cI \cap \duals)
\subseteq \cI^\sharp$. Thus in order to prove that $\cI^\sharp$ is dense in
$M_*$, it is by the previous lemma enough to prove that $(\dual)^2$ is dense in
$\dual$. But we have for all $v \in H$ with $\|v\| = 1$, $w_1,w_2 \in H$  and
$x \in M$ that $$\langle \de(x) \, W^*(v \ot w_1) , W^* (v \ot w_2) \rangle =
\langle x w_1 , w_2 \rangle  \ ,$$ which easily implies that $(\dual)^2$ is
dense in $\dual$. Hence $\cI^\sharp$ is dense in $M_*$.

Since $I \cap \duals$ is dense in $M_*$, $1$ belongs to the $\si$-strong$^*$
closure of $\lambda(I \cap \duals)^*$. Combining this with the fact that
$\lambda(I \cap \duals)$ is a $\si$-strong$^*$--norm core for $\lah$ and the
inclusion $\lambda(\cI \cap \duals)^*\lambda(\cI \cap \duals) \subseteq
\lambda(\cI^\sharp)$, we conclude that  $\lambda(\cI^\sharp)$ is a
$\si$-strong$^*$--norm core for $\lah$.
\end{proof}

\bigskip

Let us connect the modular objects of $\vfih$ to objects already constructed on
the level of $(M,\de)$.

\smallskip

We know that the operators $P$ and $J \sde J$ strongly commute and that
$\nabh^{it} = P^{it} J \sde^{it} J$ for $t \in \R$ (see lemma 8.8 and
proposition 8.9 of \cite{KV1}). Also notice that this implies for every $a \in
\Nfi$ that $\tau_t(a) \, \sde^{-it}$ belongs to $\Nfi$ and $\hat{\nab}^{it}
\la(a) = \la(\tau_t(a) \, \sde^{-it})$.

\smallskip

Put $\hat{T} = \hat{J} \, \hat{\nab}^{\frac{1}{2}}$. So $\hat{\la}(\cN_{\vfih}
\cap \cN_{\vfih}^*)$ is a core for $\hat{T}$ and $\hat{T} \hat{\la}(x) =
\hat{\la}(x^*)$ for all $x \in \cN_{\vfih} \cap \cN_{\vfih}^*$.

\medskip

\begin{lemma}
The set $\lah(\lambda(\cI^\sharp))$ is a core for $\hat{T}$.
\end{lemma}
\begin{proof}
Since $\hat{T} = \hat{J} \hat{\nab}^{\frac{1}{2}}$, the definition of $\hat{T}$
gives clearly that $\lah(\lambda(\cI^\sharp)) \subseteq
D(\hat{\nab}^{\frac{1}{2}})$. From proposition \ref{dual.prop20}, we know that
$\lah(\lambda(\cI^\sharp))$ is a dense subspace in $H$.

We now use the notation $\rho_t$ as it was introduced in notation 8.7 of
\cite{KV1}. For every $\om \in M_*$ we denote by $\rho_t(\om)$ the element in
$M_*$ defined by $\rho_t(\om)(x)=\om(\sde^{-it} \tau_{-t}(x))$. Then
$\rho_t(\cI)=\cI$ and $\xi(\rho_t(\om))=\nabh^{it} \xi(\om)$ for all $\om \in
\cI$ and $t \in \R$. If $\om \in \duals$ and $t \in \R$, it is not so difficult
to check that $\rho_t(\om) \in \duals $ and $\rho_t(\om)^* = \rho_t(\om^*)$. It
follows that $\rho_t(\cI^\sharp) = \cI^\sharp$ for all $t \in \R$, hence
$\sih_t(\lambda(\cI^\sharp)) = \lambda(\rho_t(\cI^\sharp)) =
\lambda(\cI^\sharp)$ for all $t \in \R$.

Therefore $\nabh^{it} \lah(\lambda(\cI^\sharp)) = \lah(\lambda(\cI^\sharp))$
for all $t \in \R$. We conclude from all this that $\lah(\lambda(\cI^\sharp))$
is a core for $\hat{\nab}^{\frac{1}{2}}$ (see e.g. corollary 1.21 of
\cite{JK3}) and the lemma follows.
\end{proof}

\medskip

\begin{proposition}
Consider $x \in \Nfi \cap D(S^{-1})$ such that $S^{-1}(x)^* \in \Nfi$. Then
$\la(x) \in D(\hat{T}^*)$ and  $\hat{T}^* \la(x) = \la(S^{-1}(x)^*)$.
\end{proposition}
\begin{proof}
Choose $\theta \in \cI^\sharp$. Then
\begin{eqnarray*}
& & \langle \hat{T} \lah(\lambda(\theta)) , \la(x)  \rangle = \langle
\lah(\lambda(\theta)^*) , \la(x)  \rangle = \langle \lah(\lambda(\theta^*)) ,
\la(x) \rangle = \langle  \xi(\theta^*) , \la(x) \rangle  \ .
\end{eqnarray*}
Therefore the definition of $\xi(\theta^*)$ and $\theta^*$ imply that
\begin{eqnarray*}
& & \langle \hat{T} \lah(\lambda(\theta)) , \la(x)  \rangle = \theta^*(x^*) =
\bar{\theta}(S(x^*)) = \overline{\theta(S^{-1}(x))} = \overline{\langle
\xi(\theta) , \la(S^{-1}(x)^*) \rangle}
\\ & & \spat = \langle \la(S^{-1}(x)^*)
, \xi(\theta) \rangle = \langle \la(S^{-1}(x)^*) , \lah(\lambda(\theta))
\rangle \ .
\end{eqnarray*}
Thus the previous lemma implies that $\la(x)$ belongs to $D(\hat{T}^*)$ and
$\hat{T}^* \la(x) = \la(S^{-1}(x)^*)$.
\end{proof}

\bigskip

This proposition allows us to establish easily a connection between $G$ and
$\hat{T}$. Recall that the operators $G$, $N$ and $I$ were introduced in
proposition \ref{toep.antipode} and notation \ref{notatieN}.

\begin{corollary}
We have that $\hat{T}^* = G$, $\hat{\nab} = N^{-1}$ and $\hat{J} = I$.
\end{corollary}
\begin{proof}
Using proposition \ref{toep.antipode} and the strong left invariance of $\psi$
(see proposition 5.24 of \cite{KV1}), the previous result implies easily $G
\subseteq \hat{T}^*$.

Define the subspace $C$ of $D(G)$ as $C = \langle \, \la\bigl((\psi \ot
\io)(\de(y^*)(x \ot 1))\bigr) \mid x,y \in \Nfi^* \Nps\,\rangle$.

\smallskip

Let $t \in \R$.  Remember that $\hat{\nab}^{it} \la(a) = \la(\tau_t(a) \,
\sde^{-it})$ for all $a \in \Nfi$.

Choose $x,y \in \Nfi^*\,\Nps$. Then $\sde^{it}\, \tau_t(x)$ and
$\sde^{it}\,\tau_t(y)$ belong to $\Nfi^* \, \Nps$ and $$\tau_t\bigl((\psi \ot
\io)(\de(y^*)(x \ot 1))\bigl) \, \sde^{-it} = \nu^{t} \, (\psi \ot
\io)(\de((\sde^{it}\,\tau_t(y))^*)(\sde^{it}\,\tau_t(x) \ot 1)) \ . $$
Therefore the element $\hat{\nab}^{it}\la\bigl((\psi \ot \io)(\de(y^*)(x \ot
1))\bigr) =  \la\bigl(\,\tau_t\bigl((\psi \ot \io)(\de(y^*)(x \ot 1))\bigl) \,
\sde^{-it}\,\bigr)$ belongs to $C$.

We conclude that $C$ is a dense subspace of $D(\nabh^{-\frac{1}{2}})$,
invariant under the family of operators $\nabh^{it}$ $(t \in \R)$. It follows
that $C$ is a core for $\nabh^{-\frac{1}{2}}$ and thus a core for $\hat{T}^* =
\hat{J} \nabh^{-\frac{1}{2}}$. Combining this with the fact that $G \subseteq
\hat{T}^*$, we conclude that $G = \hat{T}^*$. Now the uniqueness of the polar
decomposition implies that $\hat{\nab} = N^{-1}$ and $\hat{J} = I$.
\end{proof}

\medskip

Combining the previous corollary with proposition \ref{toep.antipode} we get
the following.

\begin{corollary}
The set $$\{\,\la(x)  \mid  x \in \Nfi \cap D(S^{-1}) \text{ such that }
S^{-1}(x)^* \in \Nfi \, \}$$ is a core for $\hat{T}^*$.
\end{corollary}

\medskip\smallskip

Recall that we introduced the GNS-construcion $(H,\io,\Ga)$ for $\psi$  by
considering $\psi$ as $\vfi_\sde$ and setting $\Ga = \la_\sde$. But $\psi$ is
by definition equal to $\vfi R$. It turns out that $\hat{J}$ connects both
pictures of $\psi$:

\begin{proposition} \label{Jhoed}
We have for all $x \in \Nps$ that $\hat{J} \,\Ga(x) = \la(R(x)^*)$.
\end{proposition}
\begin{proof}
Define the anti-unitary $U : H \rightarrow H$ such that $U \Gamma(x) =
\la(R(x)^*)$ for $x \in \Nps$. Choose $a \in \Nfi$ such that $a \in D(S^{-1})$
and $S^{-1}(a)^* \in \Nfi$.

\smallskip

For $n \in \N$, we define $e_n \in M$ such that $e_n = \frac{n}{\sqrt{\pi}}
\int \exp(-n^2 t^2) \, \sde^{it} \, dt \ $. Remember that $e_n$ is analytic
with respect to $\si$ and $\si'$, implying that $\Nfi \, e_n \subseteq \Nfi$
and $\Nps \, e_n \subseteq \Nps$

Since $\tau_s(\sde) = \sde$ we see that $\tau_s(e_n) = e_n$ for $s \in \R$,
hence $e_n \in D(\tau_{\frac{i}{2}})$ and $\tau_{\frac{i}{2}}(e_n) = e_n$. By
assumption, $a \in D(\tau_{\frac{i}{2}})$, so $a \, e_n \in
D(\tau_{\frac{i}{2}})$ and $\tau_{\frac{i}{2}}(a \,e_n) =
\tau_{\frac{i}{2}}(a)\, e_n$. Hence $\tau_{\frac{i}{2}}(a
\,e_n)\sde^{\frac{1}{2}}$ is a bounded operator and its closure equals
$\tau_{\frac{i}{2}}(a)\, (\sde^{\frac{1}{2}} e_n)$.

Define the strongly continuous one-parameter group $\kappa$ of isometries of
$M$ such that $\kappa_t(x) = \tau_t(x) \, \sde^{-it}$ for $x \in M$ and $t \in
\R$. The discussion above implies (see e.g. proposition 4.9 of \cite{JK3}) that
$a \, e_n \in D(\kappa_{\frac{i}{2}})$ and $\kappa_{\frac{i}{2}}(a \, e_n) =
\tau_{\frac{i}{2}}(a)\, (\sde^{\frac{1}{2}} e_n)$.

\smallskip

By assumption $R(\tau_{\frac{i}{2}}(a))^* = S^{-1}(a)^* \in \Nfi$, implying
that $\tau_{\frac{i}{2}}(a)$ belongs to $\Nps$. So we see that
$\kappa_{\frac{i}{2}}(a \, e_n) \sde^{-\frac{1}{2}}$ is a bounded operator and
that its closure equals $\tau_{\frac{i}{2}}(a)\,e_n \in \Npsi$. Since $\la =
\Ga_{\sde^{-1}}$, this implies that $\kappa_{\frac{i}{2}}(a \, e_n) \in \Nfi$
and $$\la(\kappa_{\frac{i}{2}}(a \, e_n)) = \Ga(\kappa_{\frac{i}{2}}(a \,
e_n)\, \sde^{-\frac{1}{2}}) = \Ga(\tau_{\frac{i}{2}}(a)\,e_n) \ .$$

We know  that we have for every $x \in \Nfi$ that $\kappa_t(x) \in \Nfi$ and
$\la(\kappa_t(x)) = \hat{\nab}^{it} \la(x)$. Since $a \, e_n \in \Nfi$ and
$\kappa_{\frac{i}{2}}(a \, e_n) \in \Nfi$, we conclude (see e.g. proposition
4.4 of \cite{JK1})  that $\la(a \, e_n) \in D(\hat{\nab}^{-\frac{1}{2}})$ and
$$\hat{\nab}^{-\frac{1}{2}} \la(a\,e_n) = \la(\kappa_{\frac{i}{2}}(a \, e_n)) =
\Ga(\tau_{\frac{i}{2}}(a)\,e_n) \ .$$ Since $(\la(a\,e_n))_{n=1}^\infty$
converges to $\la(a)$ and $(\Ga(\tau_{\frac{i}{2}}(a)\,e_n))_{n=1}^\infty$
converges to $\Ga(\tau_{\frac{i}{2}}(a))$, the closedness of
$\hat{\nab}^{-\frac{1}{2}}$ implies that $\la(a) \in
D(\hat{\nab}^{-\frac{1}{2}})$ and $$\hat{\nab}^{-\frac{1}{2}} \la(a) =
\Ga(\tau_{\frac{i}{2}}(a))\ .$$ Consequently $$U\, \hat{\nab}^{-\frac{1}{2}}
\la(a) = U \Ga(\tau_{\frac{i}{2}}(a)) = \la(R(\tau_{\frac{i}{2}}(a))^*) =
\la(S^{-1}(a)^*) = \hat{T}^* \la(a) = \hat{J}\,  \hat{\nab}^{-\frac{1}{2}}
\la(a)\ .$$ Since such elements $\la(a)$ form a core for
$\hat{\nab}^{-\frac{1}{2}} = \hat{J}\,\hat{T}^*$, such elements
$\hat{\nab}^{-\frac{1}{2}} \la(a)$ form a dense subspace of $H$. Therefore
$\hat{J} =  U$ and we are done.
\end{proof}

\medskip

The equality in the next corollary is a slight adaptation of corollary
3.6.2.(iv) of \cite{Mas-Nak}.

\begin{corollary}
We have that $\hat{J}\,J = \nu^{\frac{i}{4}}\,\,J\,\hat{J}$.
\end{corollary}
\begin{proof}
As already mentioned in section \ref{sectie1}, $J' := \nu^{\frac{i}{4}}\,J$ is
the modular conjugation of $\psi$ in the GNS-construction $(H,\io,\Ga)$. Choose
$x \in \Nps \cap D(\si_{\frac{i}{2}}')$. Since $\si_{-t} \, R = R\, \si_t'$ for
all $t \in \R$, we get that $R(x) \in D(\si_{-\frac{i}{2}})$ and
$\si_{-\frac{i}{2}}(R(x)) = R(\si_{\frac{i}{2}}'(x))$. Hence $R(x)^* \in \Nfi
\cap D(\si_{\frac{i}{2}})$ and $\si_{\frac{i}{2}}(R(x)^*) =
R(\si_{\frac{i}{2}}'(x))^*$. Combining this with the previous proposition and
the definition of the modular conjugation, we get that
\begin{eqnarray*}
 \hat{J}\,J \, \Ga(x) & = & \nu^{\frac{i}{4}}\, \hat{J}\, J' \, \Ga(x) =  \nu^{\frac{i}{4}}\, \hat{J}\,\Ga(\si_{\frac{i}{2}}'(x)^*) = \nu^{\frac{i}{4}}\,
\la(R(\si_{\frac{i}{2}}'(x)^*)^*) \\ &  = &
\nu^{\frac{i}{4}}\,\la(\si_{\frac{i}{2}}(R(x)^*)^*) =
\nu^{\frac{i}{4}}\,J\,\la(R(x)^*) = \nu^{\frac{i}{4}}\,J\,\hat{J}\,\Ga(x) \ .
\end{eqnarray*}
Therefore $\hat{J}\,J = \nu^{\frac{i}{4}}\,\,J\,\hat{J}$.
\end{proof}

\medskip

Recall that we denoted by $\nabp$ the modular operator of $\psi$ in the
GNS-construction $(H,\io,\Ga)$ and by $\hat{\nabp}$ the modular operator of
$\hat{\psi}$ in the GNS-construction $(H,\io,\hat{\Ga})$.

\begin{proposition} \label{dual.prop16}
For all $s,t \in \R$ we have the following commutation relations.
\begin{align}
& \hat{\nab}^{it}\,\nab^{is} = \nu^{ist} \,\nab^{is}\,\hat{\nab}^{it}
\tekst{and} \hat{\nabp}^{it}\,\nabp^{is} =
\nu^{ist}\,\nabp^{is}\,\hat{\nabp}^{it} \label{com1}  \\ &
\hat{\nab}^{it}\,\nabp^{is}  = \nu^{ist}\,\nabp^{is}\,\hat{\nab}^{it}
\tekst{and} \nab^{is} \nabp^{it}  = \nabp^{it} \nab^{is} \label{com2}\\ &
\hat{J}\, \nab \hat{J}  = \nabp \, , \; \, J \nab J  = \nab^{-1}   \tekst{and}
J \nabp J  = \nabp^{-1} \label{com3} \\ & \hat{J} P \hat{J}  = P^{-1}
\tekst{and} \hat{J} \sde \hat{J}  = \sde^{-1} \label{com4} \\ & P^{is} \,
\nab^{it}   = \nab^{it}\,P^{is}   \tekst{and}   P^{is} \, \nabp^{it}  =
\nabp^{it}\,P^{is} \label{com5} \\ & P^{is} \, \sde^{it}  = \sde^{it} \, P^{is}
\label{com6} \\ & \nab^{is} \, \sde^{it}  = \nu^{ist} \, \sde^{it} \,\nab^{is}
  \tekst{and}  \nabp^{is} \, \sde^{it}  = \nu^{ist} \, \sde^{it} \,
\nabp^{is} \label{com7} \\ & \hat{\nab}^{is}  \, \sde^{it}  = \sde^{it} \,
\hat{\nab}^{is}   \tekst{and} \hat{\nabp}^{is}   \sde^{it}  =  \sde^{it}\,
\hat{\nabp}^{is} \label{com8}
\end{align}

All commutation relations remain true if we remove the $\ \hat{}\ $ of
$\nab$,$\nabp$ and $J$ if there is one, add a $\ \hat{}\ $ to  $\nab$,$\nabp$
and $J$ if there is not one, replace $\nu$ by $\nu^{-1}$, replace $\sde$ by
$\sdeh$ and leave $P$ unchanged.
\end{proposition}
\begin{proof}
It is easy to check that $\nabp^{it} \la(x) = \nu^{-\frac{t}{2}}
\la(\si'_t(x))$ for all $x \in \Nfi$. We already mentioned that
$\hat{\nab}^{it} \la(x) = \la(\tau_t(x) \sde^{-it})$ for all $x \in \Nfi$ and
by definition we have $\nab^{it} \la(x) = \la(\si_t(x))$ and $P^{it} \la(x) =
\nu^{\frac{t}{2}} \la(\tau_t(x))$ for all $x \in \Nfi$. Because
$\tau_t(\sde)=\sde$ for all $t \in \R$ it is easy to verify that $P^{it} \Ga(x)
= \nu^{\frac{t}{2}} \Ga(\tau_{t}(x))$ for all $x \in \Npsi$ and by definition
we have $\nabp^{it} \Ga(x) = \Ga(\si'_t(x))$ for all $x \in \Npsi$.

\medskip

Using that all three one-parametergroups $\si$, $\si'$ and $\tau$ commute and
that $\si_t(\sde^{is})= \si'_t(\sde^{is})=\nu^{ist} \sde^{is}$ and
$\tau_t(\sde^{is}) = \sde^{is}$ for all $s,t \in \R$, it is straightforward to
check the first equality in equation (\ref{com1}), equation (\ref{com2}) and
equation (\ref{com5}) by applying the operators to an element $\la(x)$ with $x
\in \Nfi$. Using proposition \ref{Jhoed} and the fact that $\si_t R = R
\si'_{-t}$ we can check the equalities $\hat{J} \nab^{it} \hat{J} =
\nabp^{-it}$ and $\hat{J} P^{it} \hat{J} = P^{it}$ on a vector $\Ga(x)$ when $x
\in \Npsi$. This gives the first equalities of equations (\ref{com3}) and
(\ref{com4}), and the rest of equation (\ref{com3}) follows from modular
theory, because $\nu^{\frac{i}{4}}J$ is the modular conjugation of $\psi$ in
the GNS-construction $(H,\io,\Ga)$. By the biduality theorem we also get $J
\hat{\nab} J = \hat{\nabp}$ and $J P J = P^{-1}$. Because for all $t \in \R$ we
have $\nabh^{it} = P^{it} J \sde^{it} J$, we get $\hat{\nabp}^{it}=P^{-it}
\sde^{-it}$. This implies that $\hat{\nabp}^{it} \la(x) = \nu^{-\frac{t}{2}}
\la(\sde^{-it} \tau_{-t}(x))$ for all $x \in \Nfi$ and we can check then the
second equality in equation (\ref{com1}) on a vector $\la(x)$ with $x \in
\Nfi$.

\medskip

Because $R(x)=\hat{J} x^* \hat{J}$ for all $x \in M$ and $R(\sde) = \sde^{-1}$,
we get the second equality of equation (\ref{com4}). Equations (\ref{com6}) and
(\ref{com7}) follow because $P^{is},\nab^{is}$ and $\nabp^{is}$ implement
respectively $\tau_s,\si_s$ and $\si'_s$ on $M$. Also $\hat{\nab}^{is}$
implements $\tau_s$ on $M$ and this gives the first equality of equation
(\ref{com8}). Because we already saw that $\hat{\nabp}^{is} = P^{-is}
\sde^{-is}$ the second equality of equation (\ref{com8}) follows immediately
from equation (\ref{com6}).

\medskip

By the biduality theorem we can indeed perform the operation stated in the
proposition, because $P^{it} \lah(y) = \nu^{-\frac{t}{2}}
\lah(\hat{\tau}_t(y))$ for all $y \in \Nfih$ and so in a sense $\hat{P}=P$.
\end{proof}

Recall that we introduced the GNS-maps $\la$, $\Ga$, $\lah$ and $\Gah$ for the
weights $\vfi$, $\psi$, $\vfih$ and $\psih$ respectively. Using this we will
define now three new multiplicative unitaries on $H \ot H$ and relate them with
our multiplicative unitary $W$ we used all the time. Recall that we already
mentioned $V$ in section \ref{sectie1}.
\begin{definition} \label{unitaries}
Applying the von Neumann algebraic counterpart of theorem 3.16 in \cite{KV1}
(and its right invariant version) to $(M,\de)$ and $(\hat{M},\deh)$, one can
define the unitaries $V$, $\hat{W}$ and $\hat{V}$ on $H \ot H$ by the following
formulas.
\begin{align*}
V(\Ga(x) \ot \Ga(y)) &= (\Ga \ot \Ga)(\de(x)(1 \ot y)) \tekst{for all} x,y \in
\Npsi \\ \hat{W}^* (\lah(x) \ot \lah(y)) &= (\lah \ot \lah)(\deh(y)(x \ot 1))
\tekst{for all} x,y \in \Nfih \\ \hat{V} (\hat{\Ga}(x) \ot \hat{\Ga}(y)) &=
(\hat{\Ga} \ot \hat{\Ga})(\deh(x) (1 \ot y)) \tekst{for all} x,y \in
\mathcal{N}_{\psih} \ .
\end{align*}
\end{definition}
Observe that all the unitaries $W$, $V$, $\hat{W}$ and $\hat{V}$ satisfy the
pentagonal equation $$W_{12} W_{13} W_{23} = W_{23} W_{12}.$$ This could be
checked directly, but it also follows from the pentagonal equation for $W$ and
the formulas appearing in proposition \ref{unitaryW}.

\smallskip

Almost by definition we have the following.
\begin{gather*}
\de(x) = W^*(1 \ot x) W = V(x \ot 1) V^* \tekst{for all} x \in M \tekst{and} \\
\deh(y) = \hat{W}^* (1 \ot y) \hat{W} = \hat{V} (y \ot 1) \hat{V}^* \tekst{for
all} y \in \hat{M}.
\end{gather*}

The relation between all these multiplicative unitaries is given in the next
proposition.
\begin{proposition} \label{unitaryW}
We have the following formulas.
\begin{align*}
\hat{W} &= \Si W^* \Si \\ V &= (\hat{J} \ot \hat{J}) \Si W^* \Si (\hat{J} \ot
\hat{J}) \\ \hat{V} &= (J \ot J) W (J \ot J).
\end{align*}
So we have $\hat{W} \in \hat{M} \ot M$, $V \in \hat{M}' \ot M$ and $\hat{V} \in
M' \ot \hat{M}$.
\end{proposition}
\begin{proof}
The first equality follows from proposition 8.16 in \cite{KV1}. Combining lemma
8.26 in \cite{KV1} and our proposition \ref{Jhoed} we get the second equality.
Dualizing this we get $\hat{V} = (J \ot J) \Si \hat{W}^* \Si (J \ot J)$ and
this gives the third equality after applying the first one.

\smallskip

The final statement follows from the fact that $W \in M \ot \hat{M}$, the
previous formulas and the equalities $JMJ=M'$, $\hat{J} \hat{M} \hat{J} =
\hat{M}'$, $\hat{J} M \hat{J} = M$ and $J \hat{M} J = \hat{M}$. The first two
of these equalities follow from modular theory and the last two from
proposition \ref{implement}.
\end{proof}

\bigskip

\sectie{A stronger form of left invariance.}

\medskip

In this section we want to prove some stronger form of left invariance of the
Haar weight $\vfi$. We want to show that $(\io \ot \io \ot \vfi)(\io \ot
\de)(X) = (\io \ot \vfi)(X) \ot 1$ for any positive element $X \in N \ot M$ and
any von Neumann algebra $N$. The same formula is stated in \cite{EScros} for
Kac algebras, but not proved. The first proof for this formula in the Kac
algebra case was given bij Zsid\'o in \cite{Zsido} (see also remark 18.23 in
\cite{Stra}). Unfortunately the proof of Zsid\'o does not work in the case of
an arbitrary von Neumann algebraic quantum group, where possibly $\tau_t \neq
\io$.

\smallskip

In our definition of a von Neumann algebraic quantum group we assumed the
existence of invariant weights. The notion of left invariance we use, is in
fact the weakest form of left invariance that one can assume, namely
$\vfi\bigl( (\om \ot \io)\de(x) \bigr) = \vfi(x) \om(1)$ for all $\om \in
M_*^+$ and $x \in \Mfi^+$. As a special case of the next proposition we will
get the strongest form of left invariance, namely $(\io \ot \vfi)\de(x) =
\vfi(x)1$ for all $x \in M^+$. Some result in between was already proved in
proposition 5.15 of \cite{KV1}, and will be used in the proof of the
proposition.

\smallskip

When $N$ is a von Neumann algebra we denote by $N^+\Ext$ the extended positive
part of $N$ as was already mentioned in the introduction. In the proof of the
next proposition we denote by  $\langle \cdot , \cdot \rangle$ the composition
of elements in $N^+\Ext$ and $N_*^+$.

\begin{proposition}
Let $N$ be a von Neumann algebra and $X \in (N \ot M)^+$. Then we have $$(\io
\ot \io \ot \vfi)(\io \ot \de)(X) = (\io \ot \vfi)(X) \ot 1.$$ Here both sides
of the equation make sense in $(N \ot M)^+\Ext$. In particular we get $$(\io
\ot \vfi) \de(x) = \vfi(x)1$$ for all $x \in M^+$.
\end{proposition}
\begin{proof}
We will prove the proposition for the dual von Neumann algebraic quantum group
$(\hat{M},\deh)$. Because of the biduality theorem this proves the stated
result. We also represent $N$ on a Hilbert space $K$ and then it is enough to
prove the proposition in case $N = B(K)$.

\smallskip

Recall that we introduced the multiplicative unitary $\hat{V}$ in definition
\ref{unitaries}. Then define for every $z \in B(K \ot H)^+$ the element $T(z)
\in B(K \ot H)^+\Ext$ by the following formula, which makes sense because
$\hat{V} \in B(H) \ot \hat{M}$. $$T(z) = (\io \ot \io \ot \vfih) \bigl( (1 \ot
\hat{V})(z \ot 1)(1 \ot \hat{V}^*) \bigr).$$ When $\eta \in K \ot H$ we denote
by $P_\eta$ the positive rank one operator defined by $P_\eta(\xi) = \langle
\xi,\eta \rangle \eta$. Let now $\eta \in K \ot H$ and suppose $\|\eta\|=1$.
Choose an orthonormal basis $(e_i)_{i \in I}$ of $K \ot H$ such that $\eta =
e_i$ for some $i \in I$.

\smallskip

Choose $\xi \in K \ot H$. Then we have
\begin{align*}
\langle T(P_\eta) , \om_\xi \rangle & = \vfih \Bigl( (\om_\xi \ot \io) \bigl(
(1 \ot \hat{V})(P_\eta \ot 1)(1 \ot \hat{V}^*) \bigr) \Bigr) \\ & = \sum_{i \in
I} \vfih \Bigl( \bigl((\om_{\xi,e_i} \ot \io)((P_\eta \ot 1)(1 \ot \hat{V}^*))
\bigr)^* (\om_{\xi,e_i} \ot \io)((P_\eta \ot 1)(1 \ot \hat{V}^*)) \Bigr) \\ & =
\vfih \Bigl( \bigl( (\om_{\xi,\eta} \ot \io)(1 \ot \hat{V}^*) \bigr)^*
(\om_{\xi,\eta} \ot \io)(1 \ot \hat{V}^*) \Bigr).
\end{align*}
In proposition \ref{unitaryW} we saw that $\hat{V}^* = (J \ot J) W^* (J \ot J)
= (w^* \ot 1) W (w \ot 1)$ where $w = \hat{J} J$. The last equality follows
from corollary \ref{JenW}. So it follows that $$\langle T(P_\eta) , \om_\xi
\rangle = \vfih \Bigl( \bigl((\om_{(1 \ot w)\xi,(1 \ot w)\eta} \ot \io)(1 \ot
W) \bigr)^* (\om_{(1 \ot w)\xi,(1 \ot w)\eta} \ot \io)(1 \ot W) \Bigr).$$

\smallskip

In remark 8.31 of \cite{KV1} we saw that for $\om \in M_*$ one has $(\om \ot
\io)(W) \in \mathcal{N}_{\vfih}$ if and only if $\om \in \cI$. So it follows
that $\langle T(P_\eta) , \om_\xi \rangle < \infty$ if and only if
\begin{equation} \label{een}
\om_{(1 \ot w)\xi,(1 \ot w)\eta}(1 \ot \cdot) \in \cI
\end{equation}
and in that case $$\langle T(P_\eta) , \om_\xi \rangle = \| \xi(\om_{(1 \ot
w)\xi,(1 \ot w)\eta}(1 \ot \cdot)) \|^2.$$

\smallskip

Suppose that $u \in M$ is a unitary. And suppose that formula~(\ref{een}) is
valid. We claim that $$\om_{(1 \ot w)(1 \ot JuJ)\xi,(1 \ot w)\eta}(1 \ot \cdot)
\in \cI$$ and $$\xi \bigl( \om_{(1 \ot w)(1 \ot JuJ)\xi,(1 \ot w)\eta}(1 \ot
\cdot)\bigr) = R(u^*) \xi \bigl( \om_{(1 \ot w)\xi,(1 \ot w)\eta}(1 \ot \cdot)
\bigr).$$ For this choose $x \in \Nfi$ and make the following computation:
\begin{align*}
\om_{(1 \ot w)(1 \ot JuJ)\xi,(1 \ot w)\eta} & (1 \ot x^*) \\ &= \langle (1 \ot
x^*) (1 \ot \hat{J} u J) \xi,(1 \ot w)\eta \rangle \\ &= \langle (1 \ot x^*) (1
\ot R(u^*)) (1 \ot w)\xi, (1 \ot w)\eta \rangle \\ &= \langle \xi \bigl(\om_{(1
\ot w)\xi,(1 \ot w)\eta}(1 \ot \cdot)\bigr), R(u) \la(x) \rangle \\ &= \langle
R(u^*) \xi \bigl(\om_{(1 \ot w)\xi,(1 \ot w)\eta}(1 \ot \cdot)\bigr), \la(x)
\rangle
\end{align*}
From this follows our claim.

\smallskip

But then we get for every $\xi \in K \ot H$ and every unitary $u \in M$ that
$$\langle T(P_\eta) , \om_\xi \rangle = \langle T(P_\eta) , \om_{(1 \ot JuJ)
\xi} \rangle.$$ From this we may conclude that $T(P_\eta) \in (B(K) \ot
M)^+\Ext$, for all $\eta \in K \ot H$. Let now $z \in B(K \ot H)^+$. Let
$(e_i)_{i \in I}$ again be an orthonormal basis for $K \ot H$. Then $$z =
\sum_{i \in I} z^{1/2} P_{e_i} z^{1/2} = \sum_{i \in I} P_{z^{1/2} e_i}.$$ By
lower semicontinuity of $T$ we can conclude that $T(z) \in (B(K) \ot M)^+\Ext$.

\smallskip

Let now $X \in \bigl( B(K) \ot \hat{M} \bigr)^+$. Then $$T(X) = (\io \ot \io
\ot \vfih)(\io \ot \deh)(X)$$ and this clearly belongs to $(B(K) \ot
\hat{M})^+\Ext$. But it also belongs to $(B(K) \ot M)^+\Ext$ by the result in
the previous paragraph. Let $$T(X) = \infty \cdot (1-e) + \int_0^\infty \lambda
de_\lambda$$ be the unique spectral decomposition of $T(X)$, considered as an
element of $B(K \ot H)^+\Ext$. Then $$e,e_\lambda \in \bigl( B(K) \ot M \bigr)
\cap \bigl( B(K) \ot \hat{M} \bigr) = B(K) \ot \C$$ because $M \cap \hat{M} =
\C$. So take $f,f_\lambda \in B(K)$ such that $e=f \ot 1$ and $e_\lambda =
f_\lambda \ot 1$. Then define the element $S \in B(K)^+\Ext$ by $$S = \infty
\cdot (1-f) + \int_0^\infty \lambda df_\lambda.$$ Then we get that
\begin{equation}\label{vglmetS}
(\io \ot \io \ot \vfih)(\io \ot \deh)(X) = S \ot 1.
\end{equation}

\smallskip

Let us now suppose first that $K = \C$. This will prove the special case stated
in the proposition. Then $X \in \hat{M}^+$ and $S$ will be a scalar. So we get
a $\lambda \in [0,+\infty]$ such that $$(\io \ot \vfih)\deh(X) = \lambda \,
1.$$ Now there are two possibilities.
\begin{itemize}
\item Either there exists a $\om \in \hat{M}_*^+$ with $\om \neq
0$ such that $(\om \ot \io)\deh(X) \in \Mfih^+$. Then $\lambda < +\infty$
because $$\lambda \, \om(1) = \vfih\bigl((\om \ot \io)\deh(X)\bigr) < \infty.$$
But then also $$\vfih\bigl((\mu \ot \io)\deh(X)\bigr) = \lambda \, \mu(1) <
\infty$$ for all $\mu \in \hat{M}_*^+$, and so $(\mu \ot \io)\deh(X) \in
\Mfih^+$ for all $\mu \in \hat{M}_*^+$. Then it follows from proposition 5.15
in \cite{KV1} that $X \in \Mfih^+$ and so $\lambda = \vfih(X)$ because of left
invariance.
\item Either we have $\vfih \bigl( (\om \ot \io) \deh(X) \bigr) =
+\infty$ for all $\om \in \hat{M}_*^+ \setminus\{0\}$. This means that
$\lambda=+\infty$. Because of left invariance we cannot have $X \in \Mfih^+$
and so $\vfih(X)=+\infty$. Again $\lambda = \vfih(X)$.
\end{itemize}
In both cases we arrive at $(\io \ot \vfih)\deh(X) = \vfih(X) \, 1$.

\smallskip

Now we return to the general case. Let $\om \in B(K)_*^+$ and $\mu \in
\hat{M}_*^+$. Then we apply $\om \ot \mu$ to equation~(\ref{vglmetS}). This
gives us
\begin{align*}
\langle S , \om \rangle \; \mu(1) = \vfih \bigl( (\om \ot \mu \ot \io)(\io \ot
\deh)(X) \bigr) &= \vfih \bigl( (\mu \ot \io) \deh( (\om \ot \io)(X) ) \bigr)
\\ &= \mu(1) \vfih \bigl( (\om \ot \io)(X) \bigr) = \mu(1) \; \langle (\io \ot
\vfih)(X) , \om \rangle.
\end{align*}
In this computation we used the special case of the proposition proved above.
So it follows that $S = (\io \ot \vfih)(X)$ and this gives what we wanted to
prove.
\end{proof}

\bigskip

\sectie{The opposite and the commutant von Neumann algebraic quantum group}

\medskip

Given a von Neumann algebraic quantum group $(M,\de)$, represented standardly
such that $(H,\io,\la)$ is a GNS-construction for the left Haar weight $\vfi$,
we can define two new von Neumann algebraic quantum groups called the opposite
von Neumann algebraic quantum group $(M,\de)\Op$ and the commutant von Neumann
algebraic quantum group $(M,\de)'$. With the notations introduced before we
give the following definition.
\begin{definition}
The underlying von Neumann algebra of the opposite von Neumann algebraic
quantum group $(M,\de)\Op$ is again $M$ and the comultiplication $\deop$ is
given by $\deop(x)=\flip \de(x)$ for all $x \in M$.

\smallskip

The underlying von Neumann algebra of the commutant von Neumann algebraic
quantum group $(M,\de)'$ is given by $M'$ and the comultiplication $\de'$ is
defined by $\de'(x) = (J \ot J)\de(JxJ)(J \ot J)$ for all $x \in M'$.
\end{definition}
It is easy to see that $(M,\de)\Op$ and $(M,\de)'$ are again von Neumann
algebraic quantum groups. We will now give canonical choices for the left
invariant weights and their GNS-construction. As a left invariant weight on
$(M,\de)\Op$ we take $\psi$, with GNS-construction $(H,\io,\ga)$. On $M'$ we
define the weight $\vfi'$ by $\vfi'(x)=\vfi(JxJ)$ for all $x \in (M')^+$. Then
$\vfi'$ is a left invariant weight on $(M,\de)'$, with GNS-construction
$(H,\io,\la')$, where $\la'(x)=J \la(JxJ)$ for all $x \in \cN_{\vfi'}$.

\smallskip

Given these GNS-constructions we can define the multiplicative unitaries $W\Op$
and $W'$ associated to  $(M,\de)\Op$ and $(M,\de)'$ and it is clear that, using
definition \ref{unitaries} and proposition \ref{unitaryW}, they are given by
$$W\Op = \Si V^* \Si \quad\mbox{and}\quad W' = (J \ot J)W (J \ot J) =
\hat{V}.$$ It is also clear that the unitary antipode $R\Op$ of $(M,\de)\Op$
equals $R$ and the unitary antipode $R'$ of $(M,\de)'$ is given by $R'(x) =
JR(JxJ)J$ for all $x \in M'$. So the canonical right invariant weights on
$(M,\de)\Op$ and $(M,\de)'$ are $\vfi$ and $\psi'$. Then the modular elements
$\sde\Op$ and $\sde'$ are given by $$\sde\Op = \sde^{-1} \quad\mbox{and}\quad
\sde' = J\sde J.$$ One also checks easily that $\tau_t\Op$ equals $\tau_{-t}$
and $\tau'_t(x) = J \tau_{-t}(JxJ) J$ for all $t \in \R$ and $x \in M'$.

\smallskip

Defining the unitary $w = \hat{J} J=\nu^{i/4} J \hat{J}$ it is easy to see that
$\Phi:M \rightarrow M': \Phi(x) = wxw^*$ gives an isomorphism between the von
Neumann algebraic quantum groups $(M,\de)$ and $(M,\de)'\Op$. To prove this we
only have to observe that $R(x) = \hat{J}x^* \hat{J}$ for all $x \in M$ and $(R
\ot R)\de(x) = \deop(R(x))$ for all $x \in M$.

\smallskip

We conclude this section with the following formulas.
\begin{proposition}
With the notations introduced above we have:
\begin{align*}
(M,\de)\Op\hoedje &= (M,\de)\hoedje' \\ (M,\de)'\hoedje &=(M,\de)\hoedje\Op \\
(M,\de)'\Op &= (M,\de)\Op'.
\end{align*}
\end{proposition}
\begin{proof}
Because $W\Op = \Si V^* \Si$, the von Neumann algebra underlying
$(M,\de)\Op\hoedje$ is given by $$\bigl\{ (\om \ot \io)(W\Op) \mid \om \in
B(H)_* \bigr\}\dpr = \bigl\{ (\io \ot \om)(V^*) \mid \om \in B(H)_* \bigr\}\dpr
= \hat{M}'.$$ The last equality follows from proposition \ref{unitaryW}.
Further we have for every $x \in \hat{M}'$ that $$ \deop\hoedje(x) = \Si W\Op(x
\ot 1) (W\Op)^* \Si = V^*(1 \ot x)V.$$ Because $V=(\Jh \ot \Jh)\Si W^* \Si (\Jh
\ot \Jh)$, this gives $$\deop\hoedje(x)=\Si(\Jh \ot \Jh)W (\Jh x \Jh \ot
1)W^*(\Jh \ot \Jh) \Si =(\Jh \ot \Jh) \deh(\Jh x \Jh) (\Jh \ot \Jh) = \deh'(x).
$$ This gives $(M,\de)\Op\hoedje = (M,\de)\hoedje'$.

\smallskip

Applying this last formula to $(M,\de)\hoedje$ and using the biduality theorem
we get $(M,\de)\hoedje\Op\hoedje = (M,\de)'$. Taking the dual and using once
again the biduality theorem this gives our second result $(M,\de)\hoedje\Op =
(M,\de)'\hoedje$.

\smallskip

To compute $(M,\de)\Op'$ we have to observe once again that the modular
conjugation $J'$ of the left invariant weight $\psi$ on $(M,\de)\Op$ is given
by $J'=\nu^{i/4}J$. Then it is clear that $(M,\de)\Op'=(M,\de)'\Op$.
\end{proof}

\end{document}